\documentclass[11pt]{article}
\usepackage{graphics}
\usepackage{epsfig}
\def\captionsize{\footnotesize}

\def\oneeighth{ {\textstyle{\frac{1}{8}}}}

\def\onehalf{ {\textstyle{\frac{1}{2}}}}
\def\onethird{ {\textstyle{\frac{1}{3}}}}
\def\twothirds{ {\textstyle{\frac{2}{3}}}}
\def\onefourth{ {\textstyle{\frac{1}{4}}}}

\def\threehalves{ {\textstyle{\frac{3}{2}}}}

\newcommand{\be}{\begin{equation}}
\newcommand{\ee}{\end{equation}}

\def\bfUpsilon{\mbox{\boldmath$\Upsilon$}}

 \def\calE{{\cal E}}
 \def\calG{{\cal G}}
 \def\calD{{\cal D}}

 \def\calS{{\cal S}}

 \def\calB{{\cal B}}

\def\Etotal{{\calE}_T}
\def\Egas{E_{gas}}

\def\Sfilm{S_{film}}

\def\onehalf{{\textstyle{ \frac{1}{2}}}}
\def\bfmathn{{\mbox{\boldmath  $n$}}}

\def\bfgamma{{\mbox{\boldmath  $\gamma$}}}
\def\bfxi{{\mbox{\boldmath  $\xi$}}}
\def\bfalpha{{\mbox{\boldmath  $\alpha$}}}
\def\bfn{{\mbox{\boldmath  $n$}}}

\def\trace{\mbox{\rm tr}}
\def\diver{\mbox{\rm div}}

\def\reals{\mbox{I$\!$R}}
\def\realR{\mbox{I$\!$R}}

\def\bfu{{\bf  u}}
\def\bfx{{\bf  x}}
\def\bfzero{{\bf  0}}
\def\bfg{{\bf  g}}
\def\bfgo{{\bf  g(\omega_0)}}

\def\bfy{{\bf  y}}
\def\bfv{{\bf  v}}
\def\bff{{\bf  f}}
\def\bfi{{\bf  i}}
\def\bfj{{\bf  j}}
\def\bfF{{\bf  F}}
\def\bfn{{\bf  n}}
\def\bft{{\bf  t}}
\def\bfb{{\bf  b}}
\def\bfk{{\bf  k}}

\def\bfC{{\bf  C}}

\def\bfA{{\bf  A}}
\def\bfX{{\bf  X}}
\def\bfB{{\bf  B}}

\def\nGores{{N}}

\newcommand{\sansM}{\mbox{\sf M}}
\newcommand{\sansS}{\mbox{\sf S}}
\newcommand{\sansT}{\mbox{\sf T}}
\newcommand{\sansU}{\mbox{\sf U}}
\newcommand{\bfG}{\mbox{\bf G}}
\newcommand{\bfH}{\mbox{\bf H}}
\newcommand{\bfI}{\mbox{\bf I}}
\newcommand{\bfS}{\mbox{\bf S}}
\newcommand{\Cof}{\mbox{\rm Cof}}

\newcommand{\adjugate}{\mbox{\rm adj}}

\newcommand{\YoungsModulus}{E}

\def\bfx{{\mbox{\boldmath  $x$}}}
\def\bfy{{\mbox{\boldmath  $y$}}}
\def\bfu{{\mbox{\boldmath  $u$}}}
\def\bfv{{\mbox{\boldmath  $v$}}}

\newcommand {\norm}[1]{\Vert #1 \Vert}

\newcommand{\dfrac}[2]{\frac{\displaystyle #1}{\displaystyle #2}}

\def\SET{\Omega \setminus \cup_{i=1}^\nGores \Sigma_i}

\def\Egas{E_P}

\def\Sfilm{S_f}

\def\inf{\mbox{\rm inf}}

\def\fstar{{f_T^*}}

\newtheorem{remark}{Remark}
\newtheorem{definition}{Definition}
\newtheorem{lemma}{Lemma}
\newtheorem{theorem}{Theorem}

\def\FSIZE{7cm}
\def\FSIZEs{5cm}
\def\FSIZEa{7cm}
\def\FSIZEb{9cm}
\begin{document}
\title{EXISTENCE THEOREMS FOR  THIN INFLATED  WRINKLED MEMBRANES
SUBJECTED TO A HYDROSTATIC PRESSURE}
\author{Frank Baginski\\
George Washington University\\
Washington, DC 20052\\
baginski@gwu.edu
 \and
 Michael Barg\\
 George Washington University\\
Washington, DC 20052\\
{mc$\_$barg@gwu.edu}
 \and
 William Collier\\
 5912 Seventeenth Street NW,\\
Washington, DC 20011\\
{wcollier@csc.com}}

%
%
%
\maketitle
\begin{abstract}
  In this paper, we
 establish
  rigorous existence theorems  for a mathematical model
of a thin  inflated wrinkled membrane that is subjected to a shape
dependent hydrostatic pressure load. We are motivated by the problem
of determining the equilibrium shape of a strained high altitude
large scientific balloon. This problem has a number of unique
features. The balloon is very thin (20-30~$\mu$m), especially when
compared with its diameter (over 100~meters). Unlike a standard
membrane, the balloon is unable to support compressive stresses and
will wrinkle or form folds of excess material. Our approach can be
adapted to a wide variety of inflatable membranes, but we will focus
on two types of high altitude balloons, a zero-pressure natural
shape balloon and a super-pressure pumpkin shaped balloon. We
outline the shape finding process  for these two classes of balloon
designs, formulate the problem of a strained balloon in an
appropriate Sobolev space setting, establish rigorous existence
theorems using direct methods in the calculus of variations, and
present numerical studies to complement our theoretical results.
\end{abstract}
%


\baselineskip 1.75pc
\parskip 0.5pc
\section{Introduction}

\label{sec:1}
 \setcounter{equation}{0}

Balloons play an important role in NASA's  current
scientific investigations,
including upper atmosphere research, high energy astrophysics,
stratospheric composition, meteorology, and astronomy.
With  the
development of the
Ultra Long Duration Balloon
(see, e.g., \cite{jonesA})
and
the possible uses  of balloons
in the exploration of planets in our solar system,
balloons will  play
an important role in  NASA's future  scientific endeavors
(see, e.g., \cite{jonesB}).
Furthermore, many of the
techniques
that have been developed for the analysis of balloons
can
be readily adapted to other
light-weight membrane structures, including
solar sails, inflatable rovers,
aerobots,
and gossamer spacecraft.

Large scientific balloons
are regularly flown by
NASA to carry out research in the stratosphere.
With a fully inflated
diameter of over 100~meters and a thickness of 20-30~microns,
a large scientific balloon does not behave like a standard
inflated membrane.
Because the balloon is so thin, it is unable to support
compressive stresses and will instead form folds of excess
material or wrinkle.
Although our theoretical results could be applied to partially inflated balloons
with the contact problem handled by constraints,
for the numerical simulations
considered in this paper,
the balloons will be fully or near fully inflated,
so wrinkling can occur, but   folds
cannot form.
Because the diameter of a typical balloon is so large
in comparison with the thickness of the  balloon skin,
the problem of a strained large scientific balloon falls
outside the realm of typical inflated membranes.
For a survey on inflated membranes, including wrinkling, see
\cite[Chapter V, Section T]{libaisimm}  and the
references therein.
To handle wrinkling,
one replaces the wrinkled region by a
smoothed out  pseudo-surface.  This can be done via
pseudo-constitutive relations as in
the work of Stein and Hedgepath \cite{SteinHedge}   where
a  variable Poisson ratio  is introduced.
We follow the approach of
Pipkin \cite{Pipkin}
where a ``relaxed'' strain energy is used
to model wrinkling (see also \cite{SteigPip}).
In addition to wrinkling,
the balloon problem is characterized by a number of other
features, including  large displacements, relatively small strains,
and a shape dependent hydrostatic pressure load.
Problems in nonlinear elasticity  have been well-studied
using a variety of methods (see, e.g.,  \cite{Ant}, \cite{libaisimm}),
including
direct methods in the calculus of variations
(see, e.g., \cite{ball}, \cite[Appendix A]{dacorogna}),
but our existence
results  for a  thin inflated  strained large scientific balloon
are the first rigorous analytical treatment of this class of membrane
problems.
Our
analytical results
affirm the use of numerical models based on
optimization, and help  explain the
efficacy of these types of schemes
(see, \cite{JAaustin}-\cite{PalmSprings}).

To date, the  workhorse of NASA's large scientific  balloon program
has been the zero-pressure natural shape balloon, a design  that
goes back to the 1950's (see \cite{umn}). In recent years, due to
scientists' demands for long duration mid-latitude balloon flights,
a design concept that has come to be known as the pumpkin balloon
has been in development. While the design shape of a balloon is a
theoretical target which the fully inflated balloon is intended to
assume, shape generating models are usually limited to weight and
pressure considerations and ignore straining in the film. However,
upon inflation, the real balloon envelope responds to the
differential pressure and strains. For these reasons, the analysis
of an inflated balloon usually is divided into two distinct parts,
{\em shape finding} and {\em stress analysis}. Shape finding is
carried out by the balloon designer. Typically, the fully inflated
balloon must maintain a payload at a constant altitude. At float
conditions, the shape of the balloon and its volume must be
determined. Archimedes' principle states that the upward force
(lift) generated by the lifting gas is equal to the total weight of
the balloon system. At equilibrium the lift is equal to the weight
of the air displaced by the lifting gas less the weight of the
lifting gas, i.e., $\mbox{Lift}= g\rho_{a} \omega_{0}- g\rho_{g}
\omega_{0}$, where $g$ is the gravititational constant, $\rho_{a}$
is the density of the atmosphere, $\rho_{g}$ is the density of the
lifting gas, and $\omega_{0}$ is the volume of the gas bubble.
$\rho_{a}$ and $\rho_{g}$ are
 functions of altitude, temperature, etc. and
are assumed
to be known. The specific buoyancy is defined to be
$b = g(\rho_{a} - \rho_g)$.
To emphasize that a quantity depends on the design conditions, we
add a subscript of `$d$'.
 For example, at float, we have $\mbox{Lift}= b_{d} \omega_{0,d}$.
We will assume that  for other altitudes corresponding to
specific buoyancy
$b$ and volume $\omega_0$,
the following relation holds
\begin{equation}
b \omega_0 = b_d \omega_{0,d}.
\label{Archimedes}
\end{equation}
Note, $b$ is decreasing as a function of altitude, so
(\ref{Archimedes})  means that volume is increasing as a function
of altitude.
For example, a balloon with   a sea level
volume of $\omega_{0}$
will expand to  approximately $300 \omega_0$ at  39~kilometers.

The shape finding processes for the two designs discussed here are
outlined in Section~2. More detail on shape finding
can
be found in
\cite{BaginskiSIAP}.
A byproduct of the shape finding process is a
gore pattern that the
manufacturer uses
to fabricate
the balloon. A gore
is a long flat tapered panel of film. The gores   are sealed edge-to-edge
in such a way that when fully inflated,
the complete balloon  will assume
a shape very similar to the one desired by the balloon designer.
The shape of the gore pattern is an
important input into the
stress analysis.

The second part of the balloon problem involves an analysis of a
pressurized elastic membrane. For efficiency, the balloon needs to
be as light as possible, yet it must be strong enough to operate
safely over its service life. A meaningful stress/strain  analysis
requires  knowledge of certain mechanical properties of the
structural elements (e.g.,  film thickness $t$, Young's modulus $E$,
and Poisson's ratio $\nu$). The membrane is assumed to be made of a
linearly elastic isotropic material, and because it is so thin, it
is unable to support compressive stresses, and will wrinkle instead.
Following the approach of Pipkin \cite{Pipkin}, we model wrinkling
by relaxation of the film strain energy density. Pipkin's approach
was
       adapted to large scientific balloons by
 Collier  (see, e.g.,
 \cite{CoRelax})
 and has been used in a number of subsequent papers
  (see, e.g., \cite{JAaustin}-\cite{BagSch}).
We assume that the balloon is in a quasi-static
steady equilibrium state and because all of the forces on the balloon are
conservative, the problem of determining an
equilibrium
shape can be formulated in terms of a variational principle.
This formulation lends
itself
to  theoretical analysis and
efficient numerical computation.
While the main results of this paper
are the existence theorems
in Section~\ref{sec:4},
we also include numerical solutions.

This paper is divided into six
sections.
In Section~\ref{sec:2}, we
provide some
historical background on mathematical models that are
used for
the shape finding process.
The most common   models
are derived
from the equations for an axisymmetric
membrane.
We also consider a cyclically symmetric
pumpkin balloon.
In Section~\ref{sec:3}, we present a mathematical model for
a strained balloon shape.
The membrane
is  modeled as  a  {\em nonlinearly
elastic membrane shell}  (see, \cite[Sec.~9.4]{ciarletIII}).
It will be important to distinguish at least
three types of
balloon configurations:
(a)~$\Omega \subset \realR^2$ -
the flat reference  configuration or {\em natural unstrained state};
(b)~$\calS_0 \subset \realR^3$ - an initial curved configuration;
 (c)~$\calS \subset\realR^3$  - a general configuration
 (not necessarily in equilibrium).
$\Omega$ and $\calS_0$ are generated by the shape finding process
and $\calS$ is generated in the analysis of a strained shape. Our
numerical solution process is robust  and $\calS_0$ need not be
close  to the final strained equilibrium state.

 Many
 non-shallow shell theories start with an initial curved reference
 configuration  in equilibrium for some load, and formulate the
 equilibrium equations in terms of an appropriate displacement field
 $\bfu$ from that state corresponding to a  load increment.
 Since
 the actual balloon is
 constructed from
 a number of  thin flat tapered panels  of film,
 it is more natural to parameterize our
 problem in terms of
 displacement
 from the flat reference configuration,
 i.e., $\bfx : \Omega \rightarrow \reals^3$.
 This is particularly important when considering
 pumpkin balloons, where {\em any} curved configuration is necessarily strained and/or
 wrinkled
 and there is no obvious
 equilibrium  state $\calS_0$ to begin the solution process.
 In our approach, we start with
 a shape $\calS_0$  (usually obtained
 from the shape finding process),
 but do not require  it to be in equilibrium.
Our film strain energy density function $W_f$ is equivalent to the
two-dimensional strain energy of  Koiter's nonlinearly elastic
membrane shell (see,  \cite[p.~450]{ciarletIII}). We will assume the
balloon is a linearly elastic isotropic material. Our finite element
model (FEM)  triangulates $\Omega$, then uses piecewise linear
elements and constant strain triangles in its implementation. As we
will see, closed subspaces of $W^{1,4}(\Omega)$ are the natural
setting in which to study the balloon problem. In order for a
mathematical model to produce meaningful results for this type of
membrane, it is important to take wrinkling into account. This can
be accomplished by  replacing  the usual strain energy density $W_f$
with $W_f^*$, the  quasiconvexification of $W_f$ (see
\cite{Pipkin}). $W_f^*$ is the largest convex function that does not
exceed $W_f$, and in our case, we are able to calculate an explicit
formula for $W_f^*$. Using the relaxed strain energy density of the
membrane, the existence of equilibrium balloon shapes follows from
direct methods in the calculus of variations (see, e.g.,~\cite{dacorogna}).

The problem of determining the strained state of a balloon is
formulated in Section~\ref{sec:3} and existence results are
established in Section~\ref{sec:4}. Variations on the model
presented in Section~3 have been used for numerical studies of
strained balloons (see, \cite{JAaustin}-\cite{BagSch}), however,
Section~\ref{sec:4} contains the  first rigorous existence results
for a  model of this type. The tendon model implemented in earlier
work (see, e.g., \cite{JAaustin}-\cite{BagSch}) treated the tendons
as linearly elastic strings whose strain energy was added directly
into the total energy of the balloon system. In the present paper,
we focus on inextensible tendons, so tendons are included through
the use of constraints. Using a very stiff tendon in the models
\cite{JAaustin}-\cite{BagSch}, we obtain results that are consistent
with the results presented here. To complement our theoretical
results, we will present numerical solutions of strained pumpkin and
strained zero-pressure natural shape balloons in
Section~\ref{sec:5}. Our theoretical results apply to asymmetric
balloon shapes as well as cyclically symmetric ones, but we will
limit our numerical studies here to  cyclically symmetric shapes.
See \cite{ASRparis} for examples of asymmetric balloon shapes.
Section~\ref{sec:6} includes concluding remarks.

\section{The shape finding process}

\label{sec:2}
\setcounter{equation}{0}

Before moving to the main results of this paper, it is important
to have a clear picture
of the
shape finding process and how it determines a gore pattern.
We focus on
the axisymmetric zero-pressure natural shape balloon and
the cyclically symmetric pumpkin balloon.

\subsection{Zero-pressure natural shape balloons}
\label{sec:2a}

\begin{figure}
\centerline{\mbox{\psfig{figure=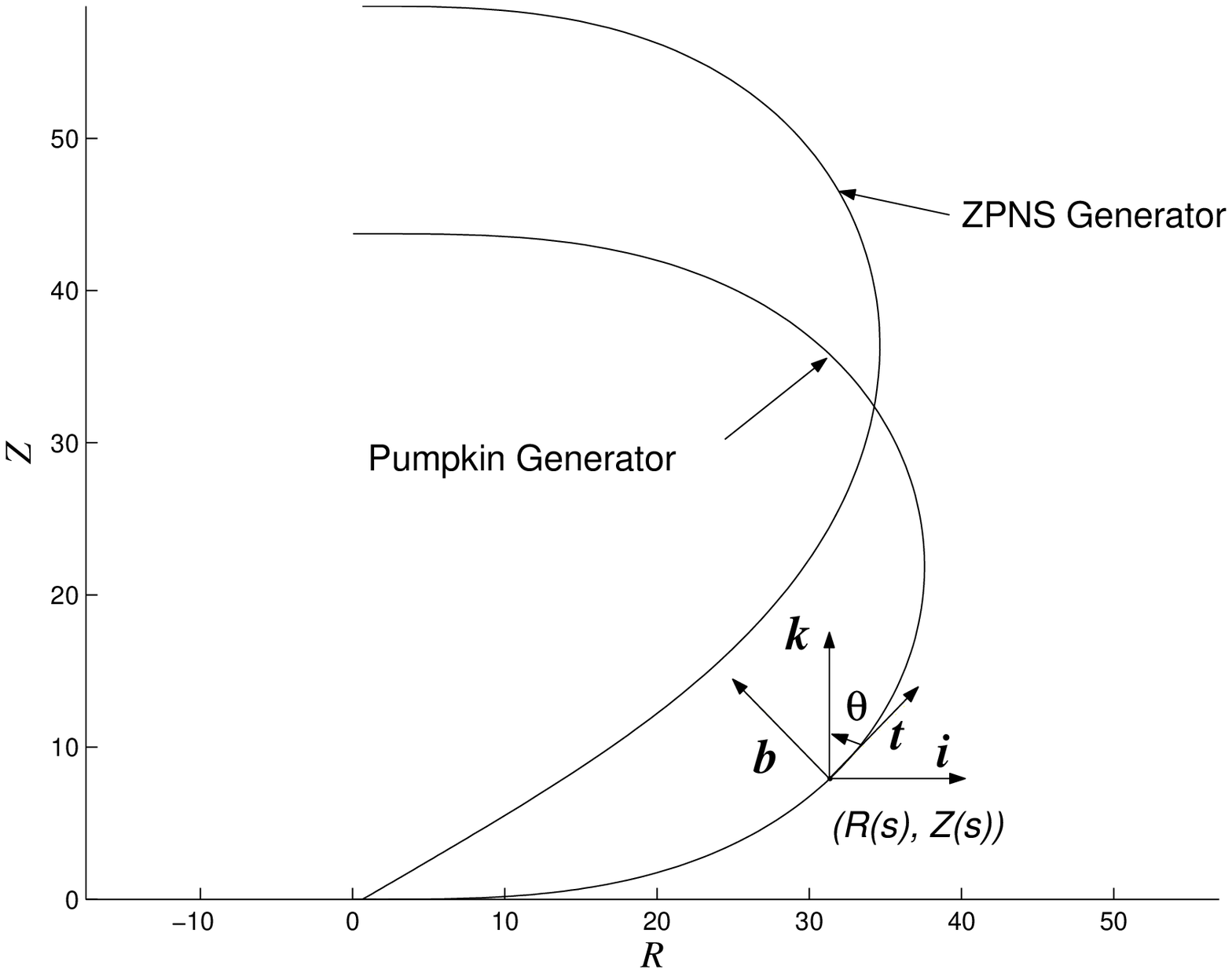,height=\FSIZE}}}
\caption{\captionsize $\bfUpsilon(s) = R(s) \bfi + Z(s) \bfk$;
(a)~Zero-pressure natural shape  gore generator;
(b)~Pumpkin gore generator  $\bfUpsilon$.
 Unit vectors
$\bft$, $\bfb$, $\bfi$, and $\bfk$  are indicated.}
\label{Upsilon}
\end{figure}

The design shape of a high altitude balloon is normally based on
conditions at  maximum (i.e., float) altitude. In addition to known
quantities such as  film weight density per unit area~$w$, suspended
payload $L$, and specific buoyancy $b$, there are other parameters
relevant to the design shape, including  the   circumferential
stress $\sigma_c$ and a constant differential pressure term $p_0$
($p_0$ is the differential pressure at the base of the balloon,
$Z=0$). The usual force  balance analysis   for an axisymmetric
membrane involving weight and differential pressure (see, e.g.,
\cite[p. 343]{Ant}) leads to a system of ordinary differential
equations
\begin{equation}
\bfzero = \frac{\partial\ }{\partial s} \left(R\sigma_m \bft \right) -
\sigma_c \bfi + R \bff
\label{ZPNSgen}
\end{equation}
for the generating curve $\bfUpsilon(s)=R(s)\bfi + Z(s) \bfk$,
where $R'(s) = \sin\theta$, $Z'(s) = \cos\theta$,
$\sigma_m$ is the meridional stress,
$\bff = -p \bfb - w \bfk$,
$p=bZ+p_0$ is the differential pressure,
$\bft$~is a unit tangent to
the generating curve, $\bfb$ is the inward normal,
and $s$ is arc length of the generating curve.
See Figure~\ref{Upsilon}(a).  $\bft$~makes
an angle of $\theta$ with
$\bfk$.
The length of the generating curve $L_d$ and $\theta(0)=\theta_0$ are unknown at the start,
and are normally found by using
a shooting method to solve a nonlinear boundary value problem
involving (\ref{ZPNSgen}) and
auxiliary conditions
(e.g.,
the total meridional tension in the $\bfk$ direction at
the bottom of the balloon must equal the weight of the suspended payload,
$L/\cos\theta_0 = 2\pi (R\sigma_m)(0)$;  zero weight at the top
of the balloon implies
  $ R(L_d)=0$).
See \cite{BaginskiSIAP} for further details.

The model that is used most commonly for applications to large scientific balloons
is the {\em zero pressure
natural shape} (ZPNS) model which assumes  $p_0=0, \sigma_c=0$
and all tension is carried in the meridional direction.
While the assumption of zero hoop stress
is clearly violated in a strained float shape,
the ZPNS design has proven to be
 successful  for NASA.
The original work on high altitude plastic balloons was
carried out in the 1950's at the University of Minnesota,
where the term natural shape balloon emerged (see, \cite{umn}).

The shape finding process for the ZPNS balloon leads to an approximation
of the axisymmetric shape by a cyclic shape where a
fundamental region $\calG_F \subset \realR^3$
is taken to be
a developable surface (see
Figure~\ref{DecoratedGore}(a)).  Thus, there is an isometry
that takes $\calG_F$
into a plane, from which  the
``lay-flat'' gore pattern $G_F$ is determined.
$G_F$ is used by the manufacturer in the construction of the balloon.
A typical large
zero-pressure natural shape balloon has over 100 gores
and it would be hard to distinguish the axisymmetric shape
from the cyclically symmetric shape generated by the $\calG_F$.
However, in
Figure~\ref{DecoratedGore}(b), we use a small number of gores
in order to highlight certain features.

\subsection{Pumpkin balloons}
\label{sec:2b}

In the past, long duration balloon
flights took place
in regions such as Antarctica where the wind currents are
regular and the balloon can
circumnavigate the South Pole under the same thermal conditions
for a long period.
In 2005,  for example,
the Cosmic Ray Energetic and Mass Experiment (CREAM)
using a ZPNS balloon
flew for nearly 42 days,
maintaining an altitude between 38 and 39 kilometers for most of its flight.
However, at
 mid-latitudes, a significant amount of ballast must be carried
 by a ZPNS balloon
 in order to  maintain  constant altitude over
 several diurnal cycles.
At night when the temperature of the lifting gas cools,  ballast
must be dropped, and  during the day when  the temperature of the
gas rises,  gas must be vented. Typically, a balloon carries enough
ballast for a few of these cycles, restricting the length of a ZPNS
mid-latitude flight to no more than a few days. One way to avoid the
need to carry significant ballast is to construct a balloon that can
contain a sufficient amount of gas to maintain altitude at night and
is strong enough to hold the overpressure caused by solar heating
during the day. This led to the consideration of a balloon design
that has come to be known as the {\em  pumpkin balloon}, a term
coined by Smalley in the early 1970's (see, \cite{Smc}). The pumpkin
shape was fore-shadowed by Taylor in his work on parachutes (see
\cite[Fig. 1]{taylor}). Currently,  NASA's Ultra Long Duration
Balloon (ULDB) Program aims to develop a pumpkin balloon that is
capable of staying aloft for one hundred days at any latitude. As of
the writing of this paper, the ULDB pumpkin has exhibited deployment
problems, but these issues are outside the scope of this paper. For
more on this topic, see \cite{JAaustin}-\cite{PalmSprings}.

\begin{figure}
 \begin{center}\captionsize
\begin{tabular}{cc}
\mbox{\psfig{figure=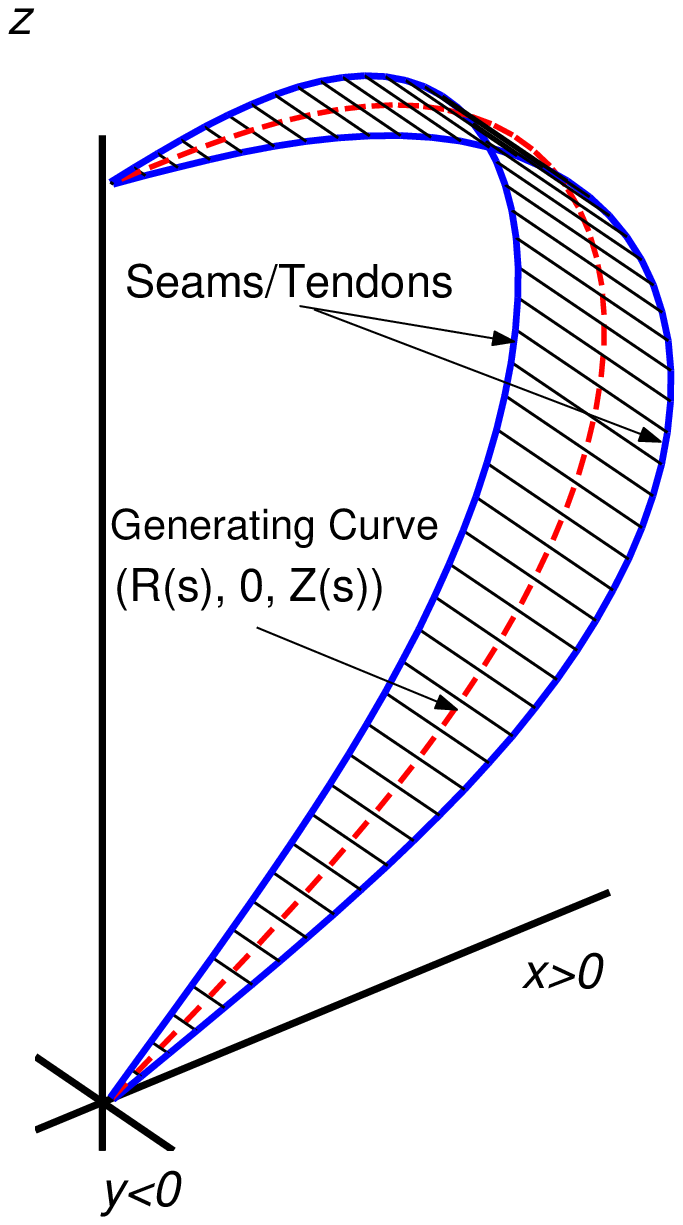,height=\FSIZE}}
&
\mbox{\psfig{figure=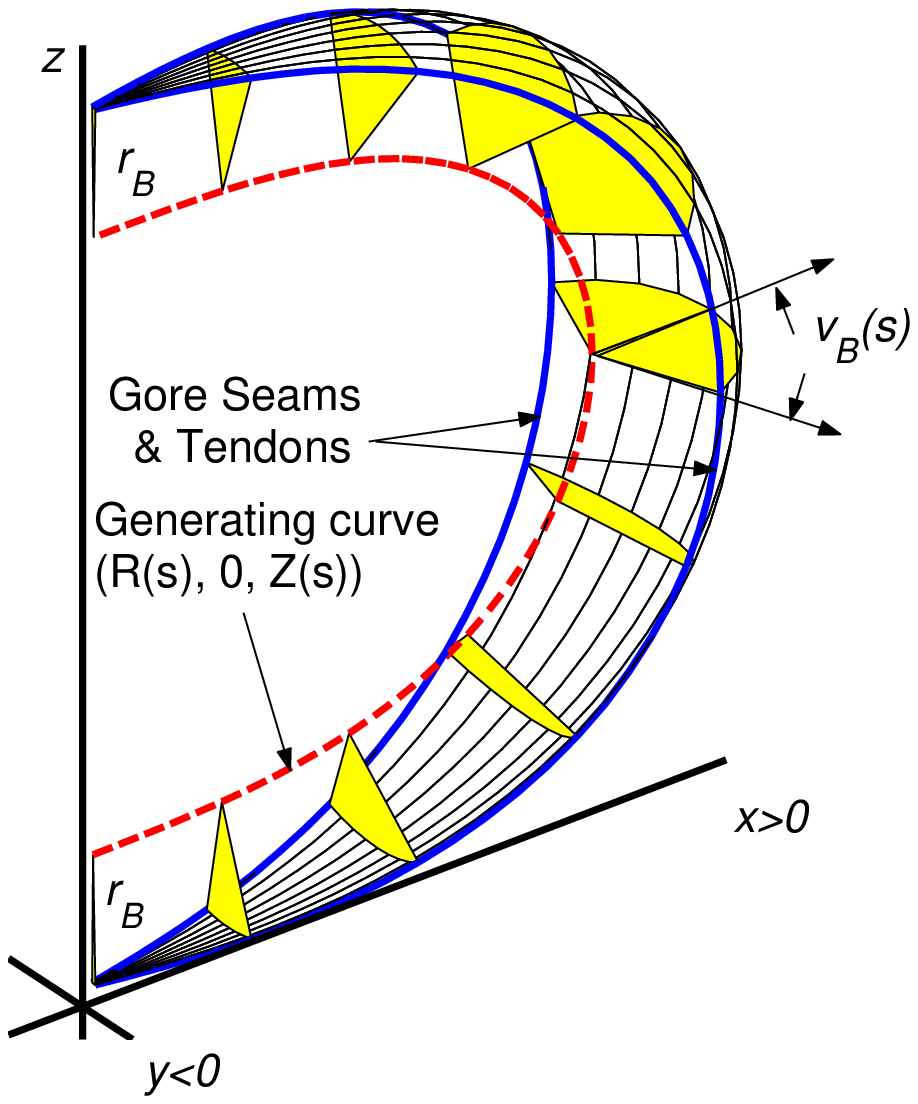,height=\FSIZE }}\cr
(a)~ZPNS~$\calG_F$   & (b)~ Pumpkin~$\calG_F$
\end{tabular}
\end{center}
\caption{\captionsize (a)~$\calG_F$ for a zero-pressure natural shape balloon;
(b)~Tubular surface $\calG_F$
with generating curve $\bfUpsilon$
and bulge radius
 $r_B$.}
\label{DecoratedGore}
\end{figure}

The principal behind the pumpkin balloon is to use a
light-weight  film as a gas barrier and  strong reinforcing
tendons
for pressure confinement and
to carry the weight of the balloon system.
Roughly speaking,
increasing the   curvature in
the circumferential direction has the effect of transferring
most of the   load to the tendons.
To model the  lobing in an ideal
pumpkin gore,  one can
assume a fundamental region $\calG_F$
is a subset of a
tubular surface. Consequently,
the
shape
$\calG_F$   is doubly curved (see Figure~\ref{DecoratedGore}).
If  a pumpkin gore
is
constructed from a flat panel of film,
any inflated pumpkin configuration
is necessarily strained and/or wrinkled.

Next, we define   pumpkin balloons that are analyzed in this paper.
We begin with a description of a constant bulge radius pumpkin that
is parameterized as a tubular surface. Let $\bfUpsilon(s)=R(s)\bfi
+ Z(s) \bfk \in \realR^3$ be a planar curve that we call the
generator of the pumpkin gore. Arc length measured along the
generating curve is denoted by $s$. See Figure~\ref{Upsilon}(b) for
a sketch of a representative $\bfUpsilon$ drawn in the $RZ$-plane.
Apriori $\bfUpsilon$ is unknown, and must be derived from
equilibrium conditions. A detailed exposition of the shape finding
equations for a pumpkin balloon, including the determination of
$\bfUpsilon$, is presented in \cite{BaginskiSIAP}. In the following,
we will assume that $\bfUpsilon$ is known. The generator is
parameterized by $s$ and a prime indicates differentiation with
respect to~$s$. Let $\bft$ denote the unit tangent  and  $\bfb =
\bft \times \bfj$ the inward unit normal of $\bfUpsilon$; $\theta =
\theta(s)$ is the angle  between $\bft$ and $\bfk$ (see
Figure~\ref{Upsilon}(b)). The set $\lbrace \bfb, \bft, \bfj
\rbrace$ gives a right-hand curvilinear basis for $\realR^3$. The
curvature of $\bfUpsilon$ is denoted by  $\kappa$ where
$\bfUpsilon''(s) = \kappa(s) \bft'(s)$. We define a tubular surface
with generator $\bfUpsilon$ and constant radius $r_B$ as follows
(see, \cite[p. 89]{doCarmo}): \be \bfx( s, v) = \bfUpsilon(s) + r_B
\left( \bfj \sin v   -\bfb(s) \cos v \right),
 \    |v| < v_B(s),  \ 0 < s < L_d,
\label{Tubular}
\ee
where
$ v_B(s)$ is
known from the shape finding process and
depends on the number of gores~$\nGores$.
See Figure~\ref{Upsilon} for a typical generating curve
and
Figure~\ref{DecoratedGore}(b)
for $\calG_F$
as determined by the shape-finding process
for a pumpkin balloon based on a tubular surface model.

\begin{figure}
\begin{center}\captionsize
\begin{tabular}{cc}
{\mbox{\psfig{figure=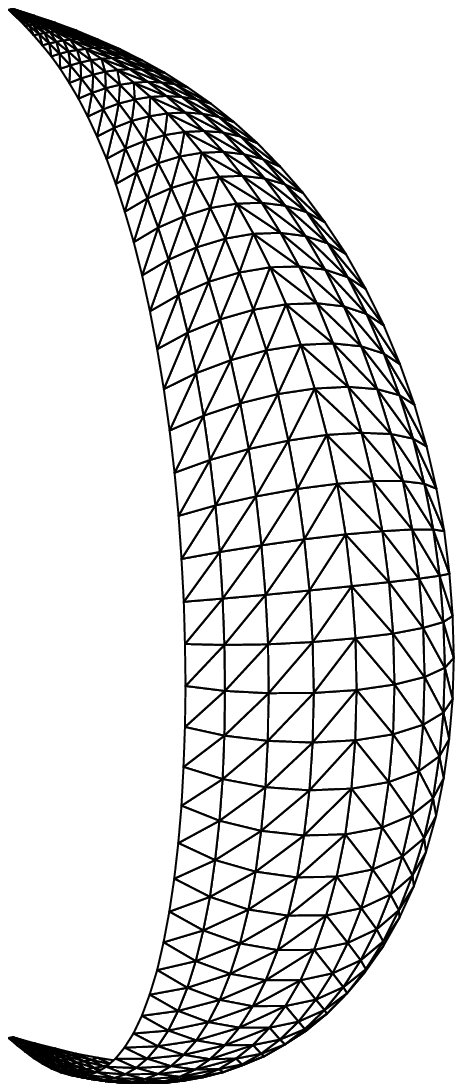,height=\FSIZEa }}
}
&
\mbox{\psfig{figure=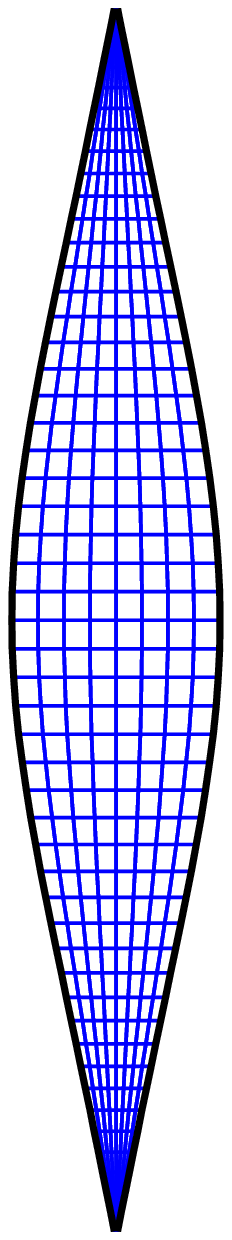,height=\FSIZEb }}\cr
{(a) $\calG_F$}   &(b)~$G_F$ -  lay-flat of gore 
\end{tabular}
\end{center}
\caption{\captionsize
(a)~$\calG_F$-initial three dimensional gore configuration;
(b)~$\Omega_F=G_F$ - usual lay flat component for one  pumpkin gore.}
\label{ThreeD}
\end{figure}
 We assume that
the pumpkin gore $\calG_F$  is
situated symmetrically with respect to the $y=0$ plane
and interior to the wedge
defined by the half-planes $y = \pm \tan(\pi/\nGores) x$
where $\nGores$ is the number of gores and
$x \ge 0$. We will refer to
$r_B$ as the {\em  bulge radius} of the pumpkin gore.
The curve traced by  $ v  \rightarrow
\bfUpsilon(s) +r_B(\bfj \sin v  -\bfb(s) \cos v  )$
is a circle lying  in the plane with normal $\bft(s)$.
By construction,
 the length
of the circular arc  that is formed
in the pumpkin gore is $2 r_B v_B(s)$
(see Figure~\ref{DecoratedGore}(b)).
We call
$v_B(s)$ the {\em bulge angle}.
For large balloons with a large number of gores,
$0< r_B\kappa(s) <<  1$, but in Figure~\ref{DecoratedGore}(b),
we used  a small number of gores in order to
 highlight
 certain features. A typical large pumpkin
balloon has
over 200 gores, yet
the lobing in a fully deployed pumpkin
is clear to see
(see \cite[Figure~2(a)]{PalmSprings}).
We define the theoretical
three-dimensional
pumpkin gore $\calG_F$ to be
\[
\calG_F = \left\lbrace (x,y,z)= \bfx(s,v),   \  -v_B(s)<   v  <
v_B(s), \
                   0 < s <  L_d
\right\rbrace.
\]
A complete shape $\calS_0$  has cyclic symmetry
and
is made of $\nGores$ copies of $\calG_F$.
Note that the length of the centerline of $\calG_F$
is
$$
{\displaystyle L_c = \int_0^{L_d} (1 + r_B \kappa(s) ) ds}.
$$
Normally, the length of a tendon is taken to be
$$
{\displaystyle L_t = \int_0^{L_d} (1 + r_B \kappa(s) \cos(v_B(s)) ds}.
$$
Corresponding to $\calG_F \subset \realR^3$ is the lay-flat configuration
$G_F \subset \realR^2$
shown in Figure~\ref{ThreeD}(b).
The respective centerlines of $\calG_F$ and $G_F$ are isometric.
The length of a
rib (i.e., a circular arc) in the
spine of $\calG_F$ is $2 r_B v_B(s)$
and this is the same as the length of a corresponding
segment orthogonal
to the centerline of $G_F$.
 It follows that $L_s$,
the corresponding
edge length
of the nominal
lay-flat pattern $G_F$,
is longer than the tendon length $L_t$. For the pumpkin example considered
in Section~5,  $ L_s = 1.005 L_t$.
A tendon is normally encased in a sleeve that is
sealed along the length of a  gore seam.
 To accommodate
the lack-of-fit,  the sleeve
is gathered first.  The tendon is tacked to the gathered
sleeve at a number of locations.
The amount of gathering is not uniform. Maximum
gathering takes place
near the equator, and decreases as one gets closer to the central axis
of the balloon.
The  smoothed out sleeve  and the edges of
two adjacent gores are heat sealed together to form one seam.
The interested
reader is referred to \cite{BagSch} and the references therein for
a more detailed discussion of tendon-foreshortening
and how it affects the film stresses.

%

\section{A model for a strained balloon}

\label{sec:3}

\setcounter{equation}{0}

In this section, we formulate
the problem of a strained balloon and
describe the components in our variational principle and related constraints.

\subsection{Preliminaries}
\label{sec:3a}
\begin{figure}
\begin{center}\captionsize
\begin{tabular}{cc}
\mbox{\psfig{figure=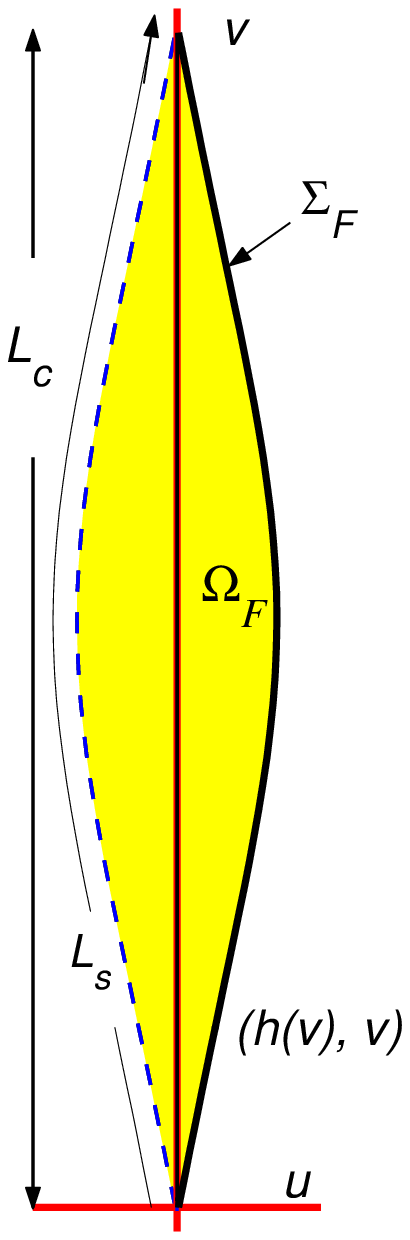,height=\FSIZE}} &
\mbox{\psfig{figure=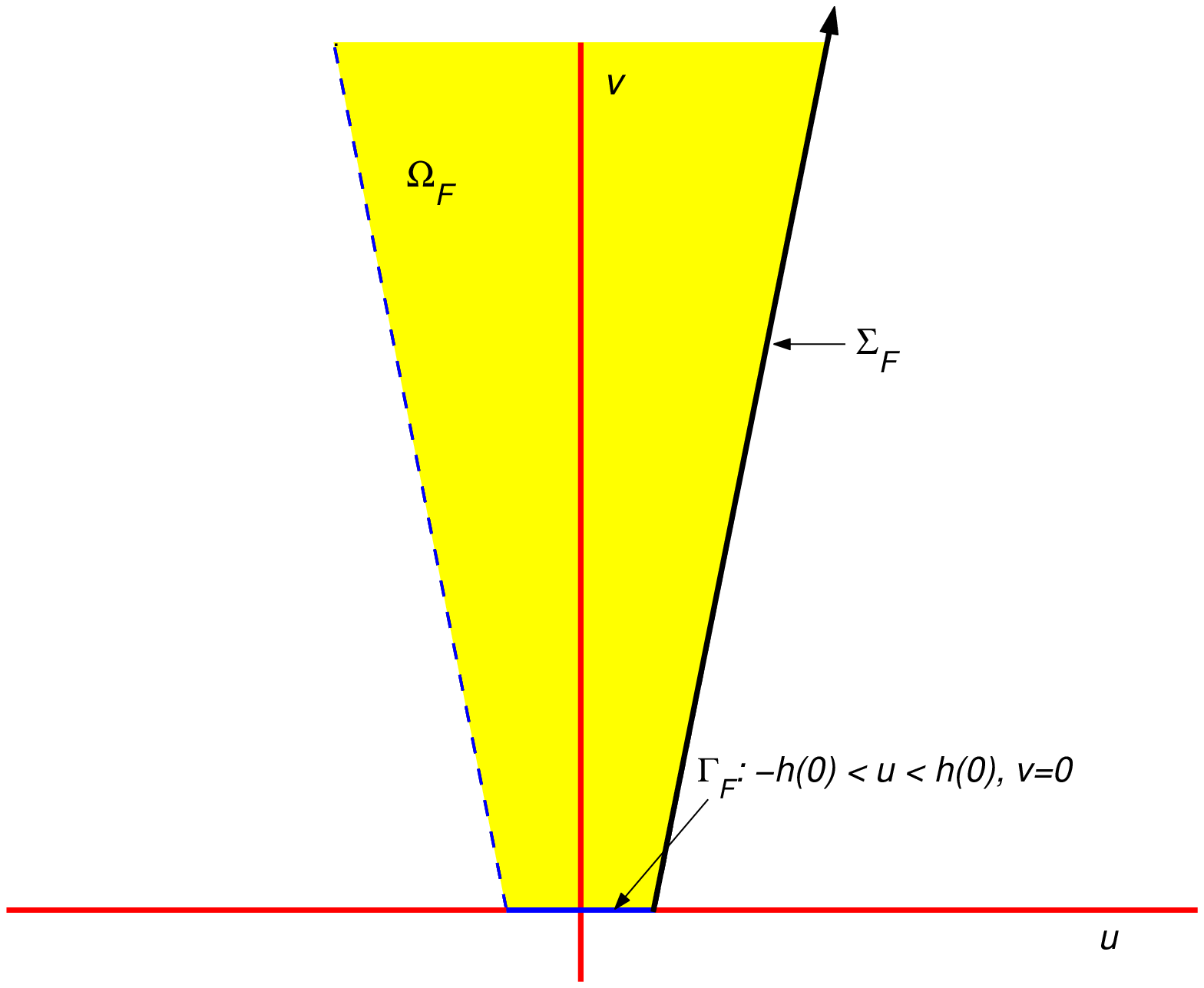,height=\FSIZE}}\cr
(a) &(b)
\end{tabular}
\end{center}
\caption{\captionsize
(a)~$\Omega_F$; seam is parametrized by $(h(v), v)$ for $ 0< v< L_c$;
(b)~close up of gore bottom;
$\Gamma_F$ component of $ \Omega_F$ that will
be  attached to the end fitting.
}
\label{OmegaF}
\end{figure}

We first define the reference configuration
for a complete balloon.
We begin with a discussion of the lay-flat configuration that is
obtained from the shape finding process.  Let
$\Omega_F=G_F$ and let $L_c$  denote the gore length where
$$
\Omega_F
= \left\lbrace (u,v) \ | \   -h(v) <  u <  h(v), \   0 <  v <  L_c
\right\rbrace,
$$
and $h(v)$ and $L_c$ are known from the shape-finding process.
We can assume that $u=h(v)$ is  a $C^{0,1}$ function, so that
$\Omega_F$ is a $C^{0,1}$ domain. Here,
$G_F$ is either a ZPNS gore or a pumpkin gore,
but in other applications it may be some other reference configuration.
In
 terms of $h(v)$,
 the seam length of $\Omega_F$ is (see Figure~\ref{OmegaF}(b))
$$
L_s = \int^{L_c}_0 {\sqrt{1 + |h'(v)|^2}}~dv.
$$
Let
\begin{eqnarray*}
\Gamma_F &=&  \left\lbrace (u,v) \ | \   -h(0) < u < h(0),  v=0
\right\rbrace.\cr
\Sigma_F &=&  \left\lbrace (u,v) \ | \    u = h(v), \ 0 <  v < L_c
\right\rbrace.
\end{eqnarray*}
Let $\triangle >0$,
$\Omega_1= \Omega_F$, $\Gamma_1=\Gamma_F$,
$\Sigma_1=\Sigma_F$, and $h_1(v) = h(v)$.
We  define
$h_i(v) = h_{i-1}(v) + \triangle$,
\begin{eqnarray*}
\Omega_i &= &
{\displaystyle \left\lbrace (u+(i-1)\triangle,v) \ | \   (u,v)\in \Omega_1
\right\rbrace} \cr
\Gamma_i &= &{\displaystyle  \left\lbrace  (u+(i-1)\triangle,v) \ | \   (u,v)\in \Gamma_1
\right\rbrace}, \cr
\Sigma_i &=&
{\displaystyle  \left\lbrace  (u + (i-1)\triangle, v)  \ | \   (u,v) \in \Sigma_1
\right\rbrace},
\end{eqnarray*}
 where   $i=2, \dots, \nGores$.
The balloon is constructed from $\nGores$ identical flat
gores $\Omega_i$
that are sealed edge to edge
(in practice, they overlap slightly  and $150\le \nGores \le 290$).
The right seam of $\Omega_i$ is identified with
the left seam of $\Omega_{i+1}$.
The right seam of $\Omega_\nGores$ is identified
with the left seam of $\Omega_{1}$. With these
identifications, we define
\be
 \Omega  = \bigcup_{i=1}^\nGores \Omega_i.
\ee

\begin{remark}
{\rm
Due to the slight overlap of the edges of
the $\Omega_i$'s as
described above, we can think of
$\Omega$ as  defining 
a two dimensional manifold. 
For convenience of exposition, we can extend slightly the gore width
of the $\Omega_i$'s  (call the extended gore $\tilde\Omega_i$)
so that adjacent local charts overlap  (i.e.,
$\tilde\Omega_i \cap \Sigma_i \cap \tilde\Omega_{i+1}  \not= \emptyset$) and $\Omega$ has
a manifold structure. 
Since the overlap is so small in our applications,
we treat it 
as a curve
$\Sigma_i$.
The metric is the usual
Euclidean metric on each $\Omega_i \subset \realR^2$.
Open sets can be defined using this metric with
special attention given to
points that lie along the seams, i.e.,
$(u,v) \in \Sigma_i \subset \Omega$.
In what follows, we study the balloon problem in
a certain Sobolev space. While it is possible to
extend the necessary definitions and results such as
the Sobolev Embedding Theorems
to a  manifold setting
(see, e.g., \cite{aubin}), introducing additional
notation for this
purpose  would be cumbersome.
Since it only makes sense to consider the differentiability
of $\bfx: \Omega \rightarrow \realR^3$ on the $\Omega_i$'s, we adhere
to the following conventions.  We write
$\bfx\in C^1(\Omega, \realR^3)$
to mean
$\bfx\in C^1(\Omega \setminus \cup_{i=1}^\nGores \Sigma_i, \realR^3)$.
With an appropriate  metric defined on $\Omega$,
we can
consider $\bfx \in C(\bar\Omega, \realR^3)$, where
the boundary of $\Omega$ is
$
\Gamma =\bigcup_{i=1}^\nGores \Gamma_i.
$
With these conventions, we  consider  mappings
$\bfx \in
C^1(\Omega, \realR^3) \cap C(\bar\Omega, \realR^3)$.}
\end{remark}

Typically, there is an end-plate located at the top and bottom of a balloon,
but
for simplicity,
we assume that a  gore comes to a point at the top.
$\Gamma$ corresponds  to the bottom end-plate.
For ease
of exposition, we will assume that the material
properties (e.g.,
Young's modulus, Poisson's ratio, film thickness, etc.)
are constant over $\Omega$.
We will formulate the equilibrium equations for a complete balloon, but
in practice, we will often assume symmetry
in order to reduce the  number of degrees of freedom in our numerical model.
We assume that
$\Omega_F =  \Omega_1$
is situated in the $uv$-plane
with its  center   along the $v$-axis
(see Figure~\ref{OmegaF}).

A deformation mapping for a  complete balloon is given by
\begin{eqnarray*}
  \bfx  : \Omega   \rightarrow  \calS \subset \realR^3,
\end{eqnarray*}
where $\calS =\bfx(\Omega)$ represents a
balloon configuration
and $\bfx \in \calD$ where
$\calD$ is the set of
admissible deformation mappings. Later,
we will
give a precise definition of $\calD$.
Let $\calS_0= \bfx_0(\Omega) $ denote the initial configuration of the deformed
balloon in
$\realR^3$.
$\calS_0$ need not be in equilibrium.
In our formulation of the problem,
the boundary of $\calS$ has one component (i.e.,
$\bfx_0(\Gamma)$,
the part that is attached to the
end-plate).
For any
$\bfx \in \calD$, we require
$$
\bfx(u,v) = \bfx_0(u,v), \ \ (u,v) \in  \Gamma.
$$
 For future reference, we note
$$
\bfalpha_i(v) = \bfx( h_i(v), v)  \mbox{\
for\ }  0<v<L_c
$$
parameterizes $\bfx|_{\Sigma_i}\subset \realR^3$,
the deformed seam
between  $\Omega_i$ and $\Omega_{i+1}$ for $i=1, \dots, \nGores-1$.
$\bfalpha_N$ parameterizes the deformed seam
between
$\Omega_\nGores$ and $\Omega_1$.

We define
\begin{eqnarray*}
| \bfx|_{1,p}^p &= & |x|_{1,p}^p + |y|^p_{1,p} + |z|^p_{1,p},
\end{eqnarray*}
where $x(u,v) = \bfx(u,v) \cdot \bfi$,  $y(u,v) = \bfx(u,v) \cdot \bfj$,
$z(u,v)= \bfx(u,v) \cdot \bfk$,
\begin{eqnarray*}
|w|^p_{1,p}   &= & \int_\Omega | w|^p~dA  + \int_\Omega |\nabla w|^p~dA,
\end{eqnarray*}
$|\nabla\bfx|^2 =  \trace\left( \nabla\bfx^T \nabla\bfx   \right) =
|\bfx_u(u,v)|^2 +  |\bfx_v(u,v)|^2$, and  $ |\bfx| = {\sqrt{ \bfx
\cdot \bfx }}$.  `$\cdot$' denotes the usual Euclidean inner product
in $\reals^3$. We will follow the convention that $\bfx_{,1}
=\bfx_{u}$ and $\bfx_{,2} =\bfx_{v}$. Let
$$
W^{1,p}(\Omega,\realR^3) = \lbrace \bfx\ : \ | \bfx|_{1,p} < \infty \rbrace.
$$

\def\SET{\Omega \setminus \cup_{i=1}^\nGores \Sigma_i}
Next, we define an  admissible deformation mapping.
\vfil\eject
 \begin{definition} Let $\bfx_0(u,v)=
 x_0(u,v) \bfi +  y_0(u,v) \bfj  + z_0(u,v) \bfk\in  W^{1,p}(\Omega,\realR^3) $.
An admissible deformation mapping,
$\bfx(u,v) = x(u,v)\bfi + y(u,v) \bfj + z(u,v) \bfk$ is a mapping
such that
\begin{description}
\item{(a)}
$\bfx\in C^1( \Omega) \cap  C(\bar\Omega)$;
\item{(b)}  $\bfx(u,v) = \bfx_0(u,v),
(u,v) \in \Gamma $;
\item{(c)}
$x(u_i,L_c) = y(u_i,L_c) = 0$, for $u_i  =(i-1)\triangle$,
\ $i=1,2,\dots,\nGores$;  \hfil\break
$z(u_i,L_c) = z(u_j, L_c)$  for all $i, j$.
\end{description}
The set  of all admissible deformation functions
is denoted $\calD$.
\label{DefA}
 \end{definition}

In the following, we let
\be
X = \left\lbrace  \bfx =
\bfx_0 +  \tilde\bfx \ | \ \tilde\bfx  \in  W^{1,p}_0(\Omega, \realR^3)
\right\rbrace.
\ee
 $X$ is  the completion of $\calD$ with respect to the norm  $| \bfx|_{1,p}$.
In our formulation of the balloon problem, $p=4$  and
a solution is
 $\bfx =   \tilde \bfx  + \bfx_0 \in X $ where
$\tilde \bfx \in W^{1,4}_0(\Omega, \realR^3)$,
and
$W^{1,4}_0(\Omega, \realR^3)$  is the completion of
$ C^1(\Omega) \cap  C_0(\bar\Omega)$ with respect to
$| \bfx|_{1,4}$.
We will  consider closed subspaces of $X$ in the form
$X_\bfg$, where $X_\bfg$ is the completion of
$$
\calD_\bfg = \
\left\lbrace
\bfx \in\calD \ | \ \bfg(u,v,\bfx, \nabla \bfx) \le \bfzero \right\rbrace
$$
with respect to $| \bfx|_{1,4}$.
$\bfg=(g_1,g_2, \dots, g_m)$   will incorporate boundary conditions and constraints
of a certain type that will ensure  $X_\bfg$ is closed.

In the following, we let
$\bfx|_\calB$ denote
the restriction of $\bfx \in \calD$ to
$\calB \subset \Omega$.
Typically, $\calB$ will denote a curve in $\Omega$ corresponding to
the middle of a gore or a gore seam.
Since $p=4$ in our application,   the trace of
$\bfx \in X$ for such a set  $\calB$ is well defined (see \cite[p. 257]{evans}).
In particular, let
$\bfx_k, \bfx\in X$ and $|\bfx_k - \bfx|_{1,4} \rightarrow 0$ and
consider $g_i$'s of the following types:
\noindent\begin{enumerate}
\item {\em Global Constraint}
\be
 g_i(\bfx) = 0, \ \ \mbox{for} \ \ g_i \in C(X, \realR).
\label{globalCOND}
\ee
From (\ref{globalCOND}),
we see ${\displaystyle \lim_{k\rightarrow \infty} g_i( \bfx_k ) =
 g_i( \lim_{k\rightarrow \infty}  \bfx_k ) =
 g_i( \bfx)}$.
\item {\em Local Constraint}
\be
g_i(\bfx|_\calB) \le 0, \ \ \mbox{for} \ \
g_i(\bfx|_\calB) \in C(\bar\calB, \realR).
\label{localCOND}
\ee
\noindent From (\ref{localCOND}), we see
${\displaystyle
\lim_{k\rightarrow \infty} g_i( \bfx_k|_{\calB} ) =
 g_i( \lim_{k\rightarrow \infty}  \bfx_k|_{\calB} ) =
 g_i( \bfx|_{\calB})
 }$.
\end{enumerate}
 We are led to the following.
\begin{lemma} Let
${\displaystyle X = \left\lbrace  \bfx =
\bfx_0 +  \tilde\bfx \ | \ \tilde\bfx  \in  W^{1,4}_0(\Omega, \realR^3)
\right\rbrace}$ and
$\bfg = (g_1, \dots, g_m)$.   If $g_i$ satisfies
 (\ref{globalCOND})  or
(\ref{localCOND}) for  $\calB= \calB_i \subset \Omega$,
$i=1,2, \dots, m$,
then $X_\bfg$ is a closed subspace, where
$$
X_\bfg
= \left\lbrace
\bfx = \bfx_0 +  \tilde\bfx \ | \ \tilde\bfx  \in  W^{1,4}_0(\Omega, \realR^3),
\ \bfg(u,v,\bfx, \nabla\bfx ) \le \bfzero \right\rbrace.
$$
\label{XgClosed}
\end{lemma}
The proof of Lemma~\ref{XgClosed}
follows directly from the
definition of $X$ and  (\ref{globalCOND})-(\ref{localCOND}).

We will demonstrate that
the components  of $\bfg$  incorporating
boundary conditions,
symmetry conditions, or
 tendon constraints satisfy (\ref{localCOND})
 and  the
 volume constraint satisfies (\ref{globalCOND}).
We first consider boundary/symmetry conditions. Tendon constraints
are considered in Section~\ref{sec:3ca} and the volume constraint is
considered in Section~\ref{sec:3cb}.

In order to
keep  the number of degrees of freedom in our numerical model
manageable,
we   use $X_\bfg$ in our simulations.
 Often, we consider
a cyclically symmetric balloon, in which case, we only need to model
one-half a gore, i.e., $y=0$ is a plane of reflectional symmetry and
the right
deformed  seam of $\Omega_1$
must lie in the  half plane
$y  = \tan(\pi/\nGores) x$ with $x\ge 0$.
Let us suppose that $\bfg$ includes boundary conditions for
a cyclically symmetric shape with $\nGores$~gores. In this case,
there are two constraint boundary conditions. Let
$\calB_1= \lbrace  (u,v) \ | \ u=0, 0 < v < L_c\rbrace$
and $\calB_2=\Sigma_1$,
\begin{eqnarray}
g_1(v,\bfx(u, v)) &= & \bfx(u,v) \cdot \bfj, \ \
(u,v) \in  \calB_1,
\label{BCb}\\
g_2( v,\bfx(u,v)) & =&\bfx( u,v ) \cdot
\left( \bfj - \tan(\pi/\nGores)  \bfi \right), \ \
(u,v) \in \calB_2. \label{BCa},
\end{eqnarray}
$g_1=0$ is the condition
that ensures
$\bfx(\Omega_1)$ is symmetric
with respect to the
$y=0$ plane
 and
$g_2=0$
is the condition
that
$\bfx(\Sigma_1)$ lies in the  plane $y  = \tan(\pi/\nGores) x $.
Suppose $\bfx_k, \bfx \in X$,
$  |\bfx_k  - \bfx|_{1,4} \rightarrow 0  $, and
 $\bfg(\bfx_k) =0$  for all $k$ where  $\bfg$ is defined in (\ref{BCb})-(\ref{BCa}).
 To demonstrate
that
$X_\bfg$ is closed,
we need to show $\bfg(\bfx) =\bfzero$, but  this is a consequence of
the fact that $g_1$ and $g_2$ are local constraints
in the form (\ref{localCOND}).
A similar argument will apply to a lobed shape with dihedral symmetry~$D_q$.
One boundary condition
is in
the form (\ref{BCb}).
If we
consider a periodic  shape with $q$ lobes and $p$ gores per lobe
where $p,q$ are positive integers and
$ p q = \nGores $, the second boundary condition
implies that $\bfx(\Sigma_q)$
lies in the $y=  \tan(\pi/q)x $ plane, i.e.,
$g_2 =\bfx(u,v) \cdot ( \bfj -\tan(\pi/q)  \bfi) =0$ for
$(u,v)\in \Sigma_q$.
In this case,  no
geometric constraints are placed
on
tendons interior to the fundamental region.

\begin{remark}
{\rm
Although
the equilibrium  state of the balloon envelope
will not have self-intersections,
we do not put this restriction
on
 $\bfx(\Omega)$ for $\bfx \in \calD$.
However,
we do require $\bfx(\Omega)$
for an equilibrium state
to be
free of self intersections.
A common  theoretical and numerical
approach  to solving a membrane problem is to begin with an
equilibrium configuration  defined by a  mapping, say
$\bfy_0$.   A small increment is then added to the
load and the corresponding displacement $\bfv^1$ is determined
for the new equilibrium, $\bfy_0+\bfv^1$. The process is continued, and
at
each stage, $\bfv^k$ is computed
for $k=1,\dots, n$.
The final solution is $\bfy_0 + \sum^n_{k=1} \bfv^k$.
Because each intermediate solution must be an equilibrium solution,
certain kinematical constraints are placed on the  $\bfv^k$'s.
In particular, each  displacement field  $\bfv^k$
cannot give rise to an
equilibrium
state with self intersections. If a solution is desired whose load is
much larger than the one corresponding to $\bfy_0$,
the number of intermediate steps could be large.
By allowing self-intersections in the class of admissible displacements,
we are able to proceed in a more direct fashion to the  desired solution.
This has
the advantage of reducing the computation time.
Although we  allow
intermediate configurations with self-intersections,
in the end, we require that the  equilibrium
configuration is free
of self-intersections.
Note, a
final solution could have regions of self-contact
(see, e.g., \cite{ASRparis}), but
this could be described via local constraints
in the form~(\ref{localCOND}).
In this paper, the balloons are
fully deployed so the contact problem does not come into play.
}
\end{remark}

\subsection{Mathematical formulation}

\label{sec:3b}

Although film and tendon weight considerations
are
relevant to
the shape finding process,
they are not significant factors in
 a stress analysis for the applications in this paper,
 and for this reason,
the film and tendon weight gravitational potential
energies are not included here.
We define
the total potential  energy of a balloon configuration  for
a
deformation mapping $\bfx$ by
 $ \Etotal = \Etotal(\bfx, \nabla\bfx)$,   where
\be
\Etotal  =\Egas +   \Sfilm,
\label{TotalEnergy}
\ee
${ \Egas}$  is the
 hydrostatic pressure potential due to the lifting gas
 and
${ \Sfilm}$  is the balloon film
strain energy.
Next, we discuss the properties of $E_P$ and $S_f$.

\subsubsection{Hydrostatic pressure potential}

\label{sec:3ba}

We follow the convention
for hydrostatic differential  pressure $P$   that
 $-P(z) >0$ means that
the inside of the balloon is pushing outward at   height $z$ units
above the base of the balloon. For a super-pressure balloon, $-P(z)
= b z + p_0$ where $b>0$ is the specific buoyancy  and $p_0>0$.
Typically,
 $p_0 >> b z_{top}$.
In a ZPNS balloon, $p_0=0$.
In our choice of coordinates, we assume
that the
base of the balloon is fixed  and corresponds to $z=0$,
so that
the potential for hydrostatic pressure  $P(z)$ is (see \cite{fisher})
\be
\Egas =  \int_{D}  P(z)  dV =
          -\int_{D} ( bz + p_0)  dV,
\label{HydrostaticPo}
\ee
where $D$ is the region occupied by the gas bubble.
Using the divergence theorem,
(\ref{HydrostaticPo}) can be replaced by
\be
\Egas =
  -  \int_\calS ( \onehalf b z^2 + p_0z) \bfk \cdot \bfn~dS,
\label{LABELEgas}
\ee
where
 $\bfn$ is the outward unit
 normal to ${ \calS}$,
 $ dS$ is surface area measure on ${\calS}$,
 and $\partial D = \calS$.
The pressure
 distribution  in an open system (i.e., without volume constraint)
 can
 be expressed as
 \begin{eqnarray}
   P( \bfx  ) & = & -(b \bfx \cdot \bfk +p_0).
   \label{PRESSURE}
 \end{eqnarray}
 Note that while $b$ and $p_0$ are known,
 the pressure potential depends on $D$.
We
 can generalize
 the discussion by allowing  a  pressure distribution  in
 the form
 $ P( \bfx )  =   \diver( \bfA(\bfx))$, but
for
(\ref{PRESSURE}), we see
\be
\bfA(\bfx) =
-\left(\onehalf b (\bfx \cdot \bfk)^2 + p_0  (\bfx \cdot \bfk)  \right)\bfk.
\label{defA}
\ee
Note, $\bfA$ is not unique.
If $P$ is given by  (\ref{PRESSURE}),
then by the divergence theorem, we have
\begin{eqnarray}
 E_P(\bfx,\nabla \bfx) &=& \int_D P(\bfx)~dV \cr
         &=& \int_{\Omega}
            f_P(\bfx, \nabla \bfx)~dA,
             \label{Pexp}
\end{eqnarray}
where
\begin{eqnarray}
f_P(\bfx, \nabla \bfx) &=& {\displaystyle g_P( \bfx, \adjugate_2 \nabla\bfx)},
\label{fPdef}\\
g_P(\bfx, \bfA ) &=&{\displaystyle -\left(\onehalf b (\bfx \cdot \ \bfk)^2 + p_0 (\bfx \cdot \ \bfk
) \right)\bfk \cdot \bfA.}
\label{gPdef}
\end{eqnarray}
By \cite[p. 117, Theorem 1.5]{dacorogna}, it follows that $f_P(\bfx,
\cdot) : \realR^{3,2} \rightarrow \realR$  is polyaffine. However,
$f_P$ polyaffine implies $f_P$ and $-f_P$ are polyconvex. Hence,
$f_P$ is quasiconvex \cite[p. 97]{dacorogna}.

Next, we establish a few inequalities.
Since
\be
|\adjugate_2\nabla\bfx \cdot \bfk|
=  | \bfx_{,1} \times \bfx_{,2} \cdot\bfk|
\leq
         | \nabla \bfx |^2,
 \label{ADJbnd}
\ee
it
follows that
\be
 -\left(\onehalf b  {|\bfx|}^2 + |p_0|~|\bfx|\right) {|\nabla \bfx|}^2  \le
    g_P(\bfx, \adjugate_2 \nabla\bfx) \le
   \left(\onehalf b {|\bfx|}^2 + |p_0|~|\bfx|\right) {|\nabla \bfx|}^2
   \label{xx1}.
\ee If $|\bfx| \le R$, then
 \be -\left(\onehalf b R^2
+|p_0|R\right)|\nabla \bfx|^2 \le g_P(\bfx, \adjugate_2 \nabla\bfx)
\le \left(\onehalf b  R^2 +|p_0|R\right) |\nabla \bfx|^2.
\label{fPBounds} \ee
\begin{remark}
{\rm
Note, for a fully inflated shape $p_0 \ge0$, but if the balloon
is not fully inflated, it is possible for $p_0 <0$.
}
\end{remark}

\subsubsection{Balloon film strain energy}

\label{sec:3bb}

Using methods in asymptotic analysis,
Ciarlet
derives the two-dimensional Koiter equations
for a nonlinearly elastic shell  in  \cite[Chapter 10]{ciarletIII}. The total energy
of a  shell of thickness $2\varepsilon$
for an appropriate displacement field $\bfu$
is
\begin{eqnarray}
\bfu & \rightarrow &
\frac{\varepsilon}{8} \int_\Omega  E^{\alpha\beta\sigma\tau,\varepsilon}
\left(
a_{\sigma\tau}(\bfu) - a_{\sigma\tau}
\right)
\left(
a_{\alpha\beta}(\bfu) - a_{\alpha\beta}
\right) \sqrt{a} dy \cr
&&
 \ \ + \ \
\frac{\varepsilon^3}{6} \int_\Omega  E^{\alpha\beta\sigma\tau,\varepsilon}
\left(
b_{\sigma\tau}(\bfu) - b_{\sigma\tau}
\right)
\left(
b_{\alpha\beta}(\bfu) - b_{\alpha\beta}
\right) \sqrt{a} dy,
\label{Koiter}
\end{eqnarray}
where $\lambda^\varepsilon$ and $\mu^\varepsilon$ are the Lam\'e
constants, and $E^{\alpha\beta\sigma\tau,\varepsilon}$ is the
two-dimensional elasticity tensor of an isotropic shell, \be
E^{\alpha\beta\sigma\tau,\varepsilon} := \frac{4 \lambda^\varepsilon
\mu^\varepsilon}{\lambda^\varepsilon + 2\mu^\varepsilon}
a^{\alpha\beta}a^{\sigma\tau} + 2 \mu^\varepsilon \left(
a^{\alpha\sigma}a^{\beta\tau} + a^{\alpha\tau}a^{\beta\sigma}
\right). \label{AMoE} \ee When applied to the balloon problem, we
can ignore the second term in (\ref{Koiter}), because $\epsilon^2$
is extremely small, i.e.,  we can ignore the  bending or flexural
energy. In our problem formulation, the natural state is the flat
reference configuration.  Thus, $a_{\alpha\beta} =
\delta_\beta^\alpha$ where $\delta_1^1= \delta^2_2=1$ and
$\delta^2_1 =\delta^1_2=0$. We can express the strain energy in
terms of a deformation from the natural state, i.e., \be \bfx
\rightarrow \frac{\varepsilon}{8} \int_\Omega
E^{\alpha\beta\sigma\tau, \varepsilon} \left( \bfx_\sigma \cdot
\bfx_\tau   -  \delta^{\sigma}_{\tau} \right) \left( \bfx_\alpha
\cdot \bfx_\beta - \delta^\alpha_\beta \right)  dA,
\label{KoiterCiarlet} \ee recognizing that $ \gamma_{\alpha\beta}=
\onehalf  \left( \bfx_\alpha \cdot \bfx_\beta  -
\delta^{\alpha}_{\beta} \right)$ is  the Cauchy-Green strain.
 Note, we do not linearize the
Cauchy-Green strain in terms of the displacement field.

The Lam\'e constants  are related to the Young's modulus $E$
and Poisson's ratio $\nu$, via
\be
\begin{array}{ccc}
\lambda &=& {\displaystyle \frac{E\nu}{(1+\nu)(1-2 \nu)}},\\
\mu     &=& {\displaystyle \frac{E}{2(1+ \nu)}},
\end{array}
\ \
\begin{array}{ccc}
E &=& {\displaystyle \frac{\mu(3\lambda + 2 \mu)}{\lambda+\mu}},\cr
\nu &=& {\displaystyle \frac{\lambda}{2(\lambda+\mu)}},
\end{array}
\label{LAM}
\ee
where we have dropped the dependence on $\varepsilon$, since we
assume that the physical constants have been determined for
a  shell of  fixed thickness.
Previous work on balloons (see, e.g., \cite{JAaustin}-\cite{BagSch})
utilized  the expression in
(\ref{KoiterCiarlet}) for the  total strain energy of
the balloon membrane.  For the convenience of the reader,
we  derive an  expression for (\ref{KoiterCiarlet})
that is consistent with the formulation in
 \cite{JAaustin}-\cite{BagSch}.

We write the
film strain energy
${ \Sfilm}$
in the form
\be
\Sfilm =
\int_\Omega W_f \thinspace dA,
\label{SfilmEQ}
\ee
where
\be
W_f =  \onehalf \bfmathn : \bfgamma.
\label{StrnEnrg}
 \ee
In (\ref{StrnEnrg}),
 ${ \bfmathn}$ represents the Second Piola-Kirchoff stress
tensor,
${\bfgamma}$ represents the  Cauchy-Green strain tensor,
and
`$:$' is the tensor inner product.
The contravariant components of ${ \bfmathn}$
are denoted by ${ n^{\alpha\beta}}$, the covariant components of
${ \bfgamma}$ are denoted by $\bfgamma_{\alpha\beta}$, and
${
\bfmathn : \bfgamma = n^{\alpha}_\beta \gamma^\beta_\alpha}.
$
Assuming a linear elastic isotropic material, we have
\be
 n^{\alpha\beta}= E^{\alpha\beta\lambda\mu} \gamma_{\lambda\mu},
 \label{LinearStressStrain}
\ee
where
${ E^{\alpha\beta\lambda\mu}}$ is the tensor of elastic moduli, i.e.,
\be
 E^{\alpha\beta\lambda\mu} =
 \frac{t \YoungsModulus}{2(1+\nu)}
\left[
a^{\alpha\lambda}\thinspace a^{\beta\mu}  +
a^{\alpha\mu}\thinspace a^{\beta\lambda}
   +
 \frac{2 \nu}{1-\nu}
a^{\alpha\beta} a^{\lambda\mu} \right], \label{MoE} \ee
$t=2\varepsilon$ is the shell thickness, and $a_{\alpha\beta} =
\delta^\alpha_\beta$. Setting $2\varepsilon=t$ and using the
relations (\ref{LAM}), we find that (\ref{AMoE}) and (\ref{MoE}) are
equivalent. In matrix form, the right Cauchy-Green deformation
tensor is
$$
\bfC = \bfF^T \bfF,
$$
where $\nabla\bfx = \bfF$ is the deformation gradient
and
the Cauchy-Green strain ($\bfgamma$) is
$$
{
\bfG = \onehalf \left( \bfC - \bfI \right).
}
$$
Assuming an isotropic film and the linear stress-strain relation
in  (\ref{LinearStressStrain}),
the
Second Piola-Kirchoff stress  can be written in matrix form  as
\be
{
\bfS =
\tau   \left[  \bfG + \nu \Cof( \bfG)^T \right],
}
\ee
where
\be
\tau  = \frac{ tE}{1-\nu^2},
\label{TAU}
\ee
and
the $2\times 2$  cofactor matrix  is
$$
\Cof\left(   \left[  \begin{array}{cc} a_{11}&a_{12}\cr
                                        a_{21}&a_{22}\end{array}
                    \right]\right) =
 \left[  \begin{array}{cc} a_{22}&-a_{12}\cr
                                       -a_{21}&a_{11}\end{array}
                    \right].
$$
$\bfG$ and $\bfS$ are symmetric
and
by the spectral representation theorem,
we have
\begin{eqnarray*}
\bfG  & = & { \delta_{1}^{\mbox{\ }} \bfn_{1}
\otimes \bfn_{1}
       + \delta_{2}^{\mbox{\ }} \bfn_{2} \otimes \bfn_{2},}   \\
{\
\bfS} & = & {\mu_{1}^{\mbox{\ }} \bfn_{1}^{\mbox{\ }} \otimes
\bfn_{1}
       + \mu_{2}^{\mbox{\ }} \bfn_{2} \otimes \bfn_{2}},
\end{eqnarray*}
where
${ \bfn_{1}}$ and $ { \bfn_{2}}$ are  orthonormal vectors.
The eigenvalues of ${ \bfS}$
 (denoted by $ \mu_{1} = \tau(\delta_1 + \nu \delta_2)$ and
 $ \mu_{2} = \tau(\delta_2 + \nu \delta_1)$)
 are  the principal stress resultants
and
the eigenvalues of ${ \bfG}$
(denoted by  ${ \delta_{1}}$ and ${ \delta_{2}}$)
 are  principal strains.
Because we have assumed a linear stress-strain constitutive
relation and an  isotropic film,
$\bfS$ and $\bfG$ have the
same principle axes.
The eigenvalues of $\bfC$ are the Cauchy strains which are
denoted by $\lambda_i^2$  where
$\delta_i  = \onehalf (\lambda_i^2 - 1)$ and
$\delta_i \ge -\onehalf$.
The film strain  density
is given by
\be
W_f =  \onehalf  \bfS : \bfG.
\label{StandardSE}
\ee
In terms of the Cauchy-Green strains $\delta_i$,
the standard membrane
strain energy is given by
 \be
 W_f(\delta_1,\delta_2)
 =
 \dfrac{tE}{2(1 - \nu^2)}
   \left({\delta_{1}}^2 + {\delta_{2}}^2 + 2 \nu \delta_{1} \delta_{2} \right),
   \label{WFb}
\ee
while in terms of the Cauchy strains,
we have
\be
  W_f(\lambda_1,\lambda_2)= \oneeighth \tau
                            \left( {\lambda_1}^4 + {\lambda_2}^4
                                    + 2\nu{\lambda_1}^2{\lambda_2}^2
             - 2(1+\nu) \left({\lambda_1}^2 + {\lambda_2}^2\right)
      + 2(1+\nu)  \right).
  \label{Wflmb}
\ee

Next, we derive a few estimates that will be needed
at a later time.
From the definition of $|\bfF|$, we have
  $$
  {|\bfF|}^2   = \lambda_1^2 + \lambda_2^2.
  $$
It follows that for $0<\nu< 1$
\be
{\lambda_1}^4 + {\lambda_2}^4
                 + 2  {\lambda_1}^2 {\lambda_2}^2
  > {\lambda_1}^4 + {\lambda_2}^4
                 + 2 \nu {\lambda_1}^2 {\lambda_2}^2
   >  \nu (\lambda_1^4 + \lambda_2^4)
                 + 2 \nu {\lambda_1}^2 {\lambda_2}^2
   =  \nu {|\bfF|}^4.\label{AAA}
\ee
From (\ref{Wflmb}) and  (\ref{AAA}), we can find a lower bound
for $W_f$, i.e.,
\begin{eqnarray}
W_f
 &\ge & {\displaystyle
 \oneeighth \tau \left( \nu |\bfF|^4 - 2(1+\nu)|\bfF|^2 + 2(1 +  \nu)
 \right).}
 \label{LOWERa}
\end{eqnarray}
Applying (\ref{AAA}) and
(\ref{LOWERa}),
we are led to
\be
\oneeighth \tau \left( \nu |\bfF|^4 - 2(1+\nu)|\bfF|^2 + 2(1 +  \nu)
 \right)
 \le W_f  \le
 \oneeighth \tau |\bfF|^4 +  \onefourth \tau(1+\nu).
 \label{WfBounds}
\ee

We conclude this subsection
with a
lemma that
will be useful
in obtaining
our existence results.
\begin{lemma}
 Let $ \kappa_1, \kappa_2 > 0$.
 \begin{description}
 \item[(i)]
 There exist
$\alpha, \rho$,  $\gamma$
such that
\be
  \begin{array}{ccccc}
    \alpha {|u|}^4 - \rho  & \leq &
     \kappa_1 {|u|}^4 - \kappa_2 {|u|}^2
     &\leq& \gamma {|u|}^4,
  \end{array}
  \label{TheoremTwo}
\ee
$0 < \alpha < \kappa_1 < \gamma$,
and $\rho =  \kappa_2^2 / 4(\kappa_1 - \alpha)$.
\item[(ii)]
There exist  constants $\gamma, \rho$ such that
\be
  \begin{array}{ccc}
     \kappa_1 {|u|}^4 + \kappa_2 {|u|}^2
     &\leq& \gamma {|u|}^4  +\rho,
  \end{array}
  \label{TheoremTwoii}
\ee
$0<  \kappa_1 < \gamma$ and
$\rho = \kappa_2^2/4(\gamma - \kappa_1)$.
\end{description}
\label{CorINEQ}
\end{lemma}
Part (i).
The second inequality in (\ref{TheoremTwo}) is obvious.
If $-\rho = \inf_u \{ (\kappa_1 - \alpha)
{|u|}^4 - \kappa_2 {|u|}^2  \} > -\infty$,
then
we are done. Otherwise,
it's easy to show that
$g(x) = (\kappa_1 - \alpha) x^4 - \kappa_2 x^2$
has an absolute minimum
at $x_0 =
\sqrt{\kappa_2 / 2(\kappa_1 - \alpha)}$ and
$g(x_0) =  - \kappa_2^2 / 4(\kappa_1 - \alpha)$. The first inequality
in (\ref{TheoremTwo}) follows with $-\rho = g(x_0)$.
The proof of Part~(ii) follows
the proof of
Part~(i).

\subsubsection{Relaxation of the film strain energy density}

\label{sec:3bc}

The energy density in  (\ref{StandardSE}) can lead to  states where
$\mu_{1}$ or $\mu_{2}$ are
negative,   corresponding  to a compression.
While this is reasonable for  certain types of shells,
the balloon film cannot support such a compression. Instead,
the film will form folds  or wrinkle.
      To tackle the problem of negative compressive stresses,
     we follow the methods introduced  by  A.~C.~Pipkin (see \cite{Pipkin}).
In the following,
let $\delta_1$ and $\delta_2$ be the principal (Cauchy-Green) strains
for a typical facet $T$ in a triangulation of $\Omega$.
 Let $\mu_1$ and $\mu_2$ be
 the corresponding principal stress resultants.
In Pipkin's approach,
a  membrane $\sansM$  is decomposed into
three distinct regions:
\begin{description}
\itemsep = 0.1 pc
\parsep = 0.0 pc
\parskip = 0.0  pc
\item{$\sansS$} - slack region ($\delta_1 < 0$,  $\delta_2<0$),
\item{$\sansT$} - tense region ($\mu_1>0$,  $\mu_2>0$), and
\item{$\sansU$} - wrinkled region ($\sansU=\sansM \setminus \sansS \cup \sansT$).
\end{description}
We apply this classification scheme
to each
$T$ in our triangulation of~$\Omega$.
On $\sansS$ the strain energy is assumed to be zero
and on $\sansT$ the relaxed strain energy density is exactly
the same as the standard strain energy  density.
On the region $\sansU$,  a modified Cauchy-Green strain
$\bfG^*$ is introduced. If $\bfG$ is the usual Cauchy-Green strain, then
\[
\bfG^* = \bfG + \beta^2 \bfn \otimes \bfn,
\]
where $\bfn$ and $\bft$ are  (unknown)
principal stress  directions
based on $\bfG^*$.
Pipkin refers to $-\beta^2 \bfn \otimes \bfn$ as the wrinkling strain
and $\bfG^*$ as the elastic strain. The elastic strain is thought to
represent the straining
in an ``averaged'' wrinkled surface and leads to uniaxial stress on
$\sansU$ in the form:
\[
\bfS^* = \mu~\bft \otimes \bft,
\]
where $\mu>0$ and $\bft$ is a unit vector orthogonal to $\bfn$.
For our exposition, we
assume that
$\bft$ is the tensile direction, and  $\bfn$ is
a unit vector orthogonal to $\bft$.
The parameter $\beta^2$ and  $\bfn$ are chosen in such a way that
the following conditions are satisfied.
\begin{eqnarray*}
\bfn \cdot \bfS^* \bfn &=&0, \\
\bfn \cdot \bfS^* \bft &=&0.
\end{eqnarray*}
For an isotropic material,
$\bfS^*$ can be written in the form:
$$
\bfS^* = \bfS +  \tau  \beta^2
\left(  \bfn \otimes \bfn + \nu \Cof( \bfn\otimes \bfn)^T \right),
$$
and  it  follows that
\be
\beta^2 = -\frac{1}{\tau} \bfn \cdot \bfS \bfn.
\label{BetaSq}
\ee
If
$W_f = W_f(\bfG)$,
its relaxation is
$W^*_f= W_f(\bfG^*)$
where $\bfG^*$ uses  $\beta^2$ from (\ref{BetaSq})
and $\bft$ is  the  principal direction
that corresponds to a positive principal strain.
One can show that
on {\sf U}, the principal strains of $\bfG^*$ are in the form
$\lbrace \delta_2, -\nu \delta_2 \rbrace$
or
$\lbrace \delta_1, -\nu \delta_1 \rbrace$
(see \cite{CoRelax}).
The principal stresses of $\bfS^* = \bfS(\bfG^*)$ are given by
\be
  \bfS^{*} = \left\{
   \begin{array}{lc}
   0,  &  {\displaystyle \delta_{1} < 0 \mbox{ and } \delta_{2} < 0},  \\
    {tE}\delta_{2} \bfn_2 \otimes \bfn_2,  &    \mu_{1} \leq 0
         \mbox{ and } \delta_{2} \geq 0, \\
  {tE}
   \delta_{1} \bfn_1 \otimes \bfn_1,  &   \mu_{2} \leq 0
         \mbox{ and } \delta_{1} \geq 0,\\
   \mu_1 \bfn_1 \otimes \bfn_1 + \mu_2 \bfn_2 \otimes \bfn_2,
    &
    \mu_{1} \geq 0 \mbox{ and } \mu_{2} \geq 0.
  \end{array}
  \right.\label{relS}
\ee
Wrinkling
is  modeled
by replacing
 $W_f$
with  $W_f^{*}$ where
\be
  W_f^{*}(\delta_1, \delta_2; t,\nu,E)  = \left\{
   \begin{array}{lc}
   0,  &  {\displaystyle \delta_{1} < 0 \mbox{ and } \delta_{2} < 0},  \\
    \onehalf{tE}
   \delta_{2}^2,  &    \mu_{1} \leq 0
         \mbox{ and } \delta_{2} \geq 0, \\
  \onehalf{tE}
   \delta_{1}^2,  &   \mu_{2} \leq 0
         \mbox{ and } \delta_{1} \geq 0,\\
   \dfrac{tE}{2(1 - \nu^2)}
   (\delta_{1}^2 + \delta_{2}^2 + 2 \nu \delta_{1} \delta_{2}),
    &
    \mu_{1} \geq 0 \mbox{ and } \mu_{2} \geq 0.
  \end{array}
  \right.\label{relWfB}
\ee

\begin{figure}
\centerline{\psfig{figure=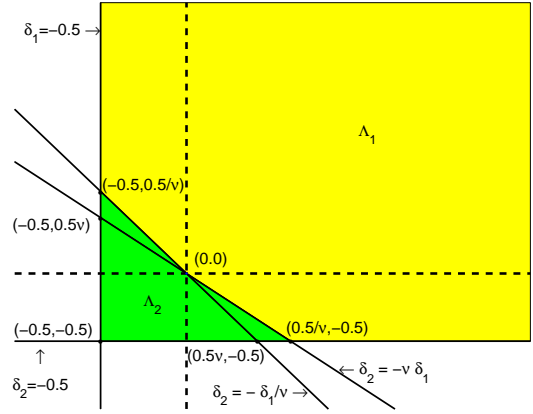,width=\FSIZE }}
\caption{\captionsize $(\delta_1,\delta_2)$ - domain for $W_f$ and $W_f^*$}
\label{LAMBDA}
\end{figure}

In order to obtain existence results, we will need upper and lower
bounds on $W_f^*$  in terms of
$|\nabla\bfx|$.
Since $W_f^*$ is the largest convex function
not exceeding $W_f$, we have
\be
W^*_f \le W_f.
\label{Upperfstar}
\ee
To obtain a lower bound for $W_f^*$,  we
consider the following.
Let
\begin{eqnarray*}
\Lambda & = & \lbrace
(\delta_1,\delta_2) \ | \ \delta_1 \ge -\onehalf, \delta_2 \ge -\onehalf
\rbrace \cr
\Lambda_1 &=& \lbrace (\delta_1, \delta_2) \in \Lambda \ | \
\delta_2 > - \delta_1/ \nu,  \delta_2 > - \nu \delta_1
\rbrace,\cr
\Lambda_2 & = & \Lambda \setminus \Lambda_1.
\end{eqnarray*}
See Figure~\ref{LAMBDA}.
While  $\Lambda_1$ is unbounded,
    $W_f =
W_f^*$ for  $(\delta_1,\delta_2)\in\Lambda_1$. Since   $\Lambda_2$ is
compact, there is a constant $d$ such that
$$
{\displaystyle
d = \max_{ (\delta_1,\delta_2) \in \Lambda_2} \lbrace
          W_f(\delta_1,\delta_2) - W_f^*(\delta_1,\delta_2)   \rbrace.
}$$
$W_f(\delta_1,\delta_2)$ has a
critical point when $\delta_1 + \nu \delta_2 =0$ and $\delta_2 + \nu \delta_1=0$.
Thus, $(\delta_1, \delta_2)=(0,0)$ is the only critical point. Hence,
$d$ is achieved on
the   boundary of $\Lambda_2$ at a corner point. By
inspection, one finds $d$ is achieved at $(-0.5,-0.5)$,  $W^*_f(-0.5,-0.5)=0$,
and $d =  \onefourth \tau (1 + \nu)$.
Thus,
\be
W_f  \le W_f^* +  \onefourth \tau (1 + \nu),  \ \  (\delta_1,\delta_2) \in \Lambda.
\label{Lowerfstar}
\ee
Combining (\ref{Upperfstar})-(\ref{Lowerfstar}), yields
\be
 W_f - \onefourth \tau (1 + \nu)
 \le W_f^*
 \le W_f.
\label{UpLowB}
\ee

We conclude this subsection by noting if
$f$ is quasi-convex,
so  is $f + C$ for any $C \in \realR$,
for if
$$\int_D f \left( u_0,v_0, \bfx_0, \bfxi_0 + \nabla\bfx(y)\right) dy \ge
|D| f( u_0,v_0, \bfx_0),
$$
for every cube $D\subset \Omega$, for every $(u_0,v_0,
\bfx_0,\bfxi_0) \in \Omega \times \realR^m \times \realR^{nm}$ and
for every $\bfx\in W^{1,\infty}_0(D;\realR^m)$, then
$$
\int_D \left( f\left( u_0,v_0, \bfx_0, \bfxi_0 + \nabla\bfx(y)\right) + C\right) dy \ge
|D| \left( f( u_0,v_0, \bfx_0) + C \right).
$$

\def\inf{\mbox{\rm inf}}

\subsection{Constraints}

\label{sec:3c}

\subsubsection {Tendon constraint}

\label{sec:3ca}

Tendons are included in our balloon model
through the use of constraints.
We assume that a
tendon
is encased in a sleeve that is attached along the seam
between two adjacent gores. We assume that the
tendons  are connected to a common
point at the top of the balloon
and to a bottom end-plate. The top of the balloon
is free to move up and down the $z$ axis and the bottom end-plate is fixed
at the level $z=0$.

\begin{remark}
{\rm
The interaction between the tendon and film is very complicated.
Even without tacking to the sleeve, a tendon will not slide freely
within its sleeve when a pumpkin balloon is fully inflated. The same
can be said for a zero-pressure balloon which uses a load tape, a
one-dimensional structural element that serves the same function as
a load tendon. For simplicity, we have used the same tendon model
for both the ZPNS and the pumpkin balloon. The assumption of
constant tension enables us to include stiff tendons in the present
model and analyze equilibrium configurations of a strained fully
inflated shape.
}
\end{remark}

We first consider a single tendon in the deformed state
$\bfx(\Sigma_i)$ for
$\bfx \in \calD$.
By the Sobolev Embedding Theorem
(see, e.g.,  \cite[Theorem 7.26, p. 171]{GiTr};
for a  manifold setting, see \cite[Theorem 2.34, p. 55]{aubin}),
 $W^{1,4}(\Omega)$ is compactly embedded in $C^{0,\beta}(\bar\Omega)$
for any $0<\beta < 1/2$.
Hence,
$\bfx \in X$ (or $X_\bfg$) implies that
$\bfx \in C^{0,\beta}(\bar\Omega, \realR^3)$,
and in particular,
$\bfalpha_i(v) =\bfx|_{\Sigma_i}$  is well-defined  and  H\"older continuous.
The length of $\bfx(\Sigma_i)$  is denoted by
$L[\Sigma_i]$.
A tendon
is under constant tension, and
the
strain in the $i$th tendon
is
\be
\epsilon_i =\left( L[\Sigma_i] - L_t \right) /  L_t.
\label{TendonStrain}
\ee
In our numerical model,
$L[\Sigma_i]$ is  approximated by
$$
L[\Sigma_i] \approx \sum_{j=1}^n |\bfalpha_i(v_j) - \bfalpha_{i}(v_{j-1})|
$$
where
$\lbrace v_j \rbrace $ is a partition of $[0, L_c]$, i.e.,
$0=v_0 < v_1 < \cdots v_{n-1} < v_n =L_c$.

The strain energy in a tendon of length $L_t$ with
stiffness $K$ is
$S_{t,i} = S_{t}[\Sigma_i] = \onehalf K {\epsilon_i}^2 L_t$.
$K$ has the units of force.  Since the tendon cannot
support compression, we utilize
the relaxed
tendon strain energy
\be
S_{t,i}^* = \left\lbrace
\begin{array}{cl}
\onehalf K {\epsilon_i}^2 L_t,  & \epsilon_i > 0,\cr
0, &\epsilon_i \le 0.
\end{array}\right.
\label{rTendonStrain}
\ee
Since we assume that the tendons are inextensible,
we  set $K=1$ and
impose
the conditions
\be
S_{t,i}^* \le 0, \ \ i=1,2,\dots, \nGores,
\label{tCON}
\ee
which can be
included as  local constraints
in the form of (\ref{localCOND}).
Because we are not using $\nabla\bfx$ to calculate the
arc length of tendons, a continuity argument similar to the
one used in Section~\ref{sec:3}
shows
$X_\bfg$ is closed when
$\bfg \le \bfzero$ includes a constraint  in the
form
(\ref{tCON}).
Consider
(\ref{TotalEnergy}) subject to
 constraints in
 the form of (\ref{tCON}).
 Let $\bfx$ be  a minimizer of this system where
   $\lambda_i, i=1, \dots, \nGores$
   are the corresponding Lagrange multipliers.
  We can interpret the $\lambda_i$'s in
\be
E_T(\bfx, \nabla\bfx) +  \sum_{i=1}^\nGores \lambda_i S_{t,i}^*.
\label{AUG}
\ee
as the tendon stiffnesses
that are needed to maintain equilibrium.
In our numerical simulations presented here,
we will consider cyclically symmetric solutions so
only one tendon constraint is needed.
Note, tendons are
modeled in a different fashion in
 \cite{JAcrystalcity}-\cite{BagSch},
but our results here are consistent with
those works when  a very large value of $K$ is  used.
We carried out numerical simulations
when (\ref{tCON}) is replaced by an equality constraint
and  obtained  essentially the same results
as the   `$\le$' case
with $\epsilon_1 \approx 10^{-10}$~m/m.
If we drop the requirement
that the tendons are inextensible,
we can replace  the `$\le$' in
 (\ref{tCON}) with  `$\ge$'. In this case,
 we find  $S_{t,i}^*>0$
 for the solution where
 $\epsilon_1 \approx 0.005$~m/m.

\subsubsection{Volume constraint}

\label{sec:3cb}

In an open system such as the ZPNS design,
the balloon is open to the atmosphere at $z=0$.
This means that the volume
of the balloon at equilibrium will adjust itself accordingly.
In a closed system,
the volume is fixed, and
$p_0$ is determined.
Using the divergence theorem,
the volume of the balloon $\calS$  can be expressed as
\be
V(\bfx )  =
\int_{D} 1~dV = \onethird \int_{\calS}  \bfx\cdot\bfn~dS
=
\int_\Omega f_V(\bfx,\nabla\bfx)~dA
\label{Vol}
\ee
where
 $\bfn~dS = \adjugate_2\nabla\bfx~dA$ and
 $f_V(\bfx, \nabla\bfx) = g_V(\bfx, \adjugate_2\nabla\bfx)$,
\be
g_V(\bfx, \bfA) =
\onethird\bfx\cdot \bfA.
\label{densV}
\ee
The same argument used
to show $f_P$ is quasiconvex, shows
$f_V(\bfx, \cdot)$ is quasiconvex.
(\ref{Vol}) can be evaluated exactly when
$\calS$ is a faceted surface.
The volume constraint is
in the form,
\be
V(\bfx) - \omega_0 = 0,
\label{VolumeConstraint}
\ee
where
$\omega_0$ is the target volume.

In a closed system where (\ref{VolumeConstraint}) holds and $p_0$ is
unknown, we can assume $-P(z) = b z$. After finding an equilibrium
solution by minimizing (\ref{TotalEnergy}) subject to
(\ref{VolumeConstraint}), we find that the Lagrange multiplier
$\lambda$ in
$$
\Etotal(\bfx, \nabla\bfx)  + \lambda (V(\bfx)- \omega_0)
$$
yields the appropriate
constant pressure term  that is needed for equilibrium.
Although we use Lagrange multipliers to handle  constraints in our
numerical approach, our theoretical model incorporates
such a  constraint  directly
into the underlying solution space.

Before we demonstrate that the volume constraint can be included in
the definition of $X_\bfg$, we will need an estimate  for
$|\adjugate_2\bfA -\adjugate_2\bfB|$.

\begin{lemma} Let $\bfA, \bfB \in \realR^{nm}$.
There exists a constant $\tilde\beta>0$ independent of $\bfA, \bfB$ such that
\be
|\adjugate_2\bfA -
\adjugate_2\bfB| \leq  \tilde \beta(|\bfA| + |\bfB|)|\bfA-\bfB|.
\label{ADJlemmaEq}
\ee
\label{ADJlemma}
When $nm=6$, $\tilde\beta = 3$.
\end{lemma}
We establish the proof for the case
$nm=6$, but
it is straightforward
to modify the argument for general $nm$.
For
$\left(
\begin{array}{cc}
x_1 & x_2\cr x_3 &x_4\cr x_5& x_6 \end{array} \right) \in
\realR^{3,2}$,  consider $\bfX =[x_1, x_2, \dots, x_6]^T$. Define
$\bff(\bfX) = [ f_1(\bfX), f_2(\bfX), f_3(\bfX)]^T = \adjugate_2\bfX
$ where $f_1(\bfX) = x_3x_6 -x_4x_5$, $f_2(\bfX) =  -(x_1x_6 -
x_2x_5)$, $f_3(\bfX) = x_1x_4 - x_2x_3$. For $f_1$, we see $D
f_1(\bfX) = [ 0, 0, x_6, -x_5, -x_4,x_3]^{T}$. Since $D^2 f_1(\bfX)$
is independent of $\bfX$, we write $D^2 f_1(\bfX) = \bfH_1$, where
$$
\bfH_1
=
\left( \begin{array}{cccccc}
    0 & 0 & 0 & 0 & 0 & 0 \\
    0 & 0 & 0 & 0 & 0 & 0 \\
    0 & 0 & 0 & 0 & 0 & 1 \\
    0 & 0 & 0 & 0 & -1 & 0 \\
    0 & 0 & 0 & -1 & 0 & 0 \\
    0 & 0 & 1 & 0 & 0 & 0
  \end{array}
\right).
$$
One easily obtains
${\displaystyle
f_1(\bfX) =
\onehalf\bfX^{T}\bfH_1\bfX}$.
For $i=2,3$, we
find
$D^2 f_i(\bfX) = \bfH_i$  and
$f_i(\bfX) = \onehalf\bfX^T\bfH_i\bfX$,
where
$$
\bfH_2
=
\left(
  \begin{array}{cccccc}
    0 & 0 & 0 & 0 & 0 & -1 \\
    0 & 0 & 0 & 0 & 1 & 0 \\
    0 & 0 & 0 & 0 & 0 & 0 \\
    0 & 0 & 0 & 0 & 0 & 0 \\
    0 & 1 & 0 & 0 & 0 & 0 \\
    -1 & 0 & 0 & 0 & 0 & 0
  \end{array}
\right)
\ \ \mbox{and} \ \
\bfH_3
=
\left(
  \begin{array}{cccccc}
    0 & 0 & 0 & 1 & 0 & 0 \\
    0 & 0 & -1 & 0 & 0 & 0\\
    0 & -1 & 0 & 0 & 0 & 0\\
    1 & 0 & 0 & 0 & 0 & 0 \\
    0 & 0 & 0 & 0 & 0 & 0 \\
    0 & 0 & 0 & 0 & 0 & 0
  \end{array}
\right).
$$

Using the
properties  of $\bff(\bfA)$ and $\bff(\bfB)$,
we obtain
\begin{eqnarray*}
|\bff(\bfA) - \bff(\bfB)| & = &\displaystyle{ \left| \left(f_1(\bfA) - f_1(\bfB),
f_2(\bfA) - f_2(\bfB), f_3(\bfA) - f_3(\bfB)\right)\right|}\\
& = & \onehalf\left| \left(\bfA^{T}\bfH_1(\bfA - \bfB) +
         (\bfA-\bfB)^{T}\bfH_1\bfB,\right.\right. \\
        & &\mbox{\hskip 0.75cm}   \bfA^{T}\bfH_2(\bfA - \bfB) +
          (\bfA-\bfB)^{T}\bfH_2 \bfB,  \\
        & &\mbox{\hskip 0.75cm}  \left. \left.\bfA^{T}\bfH_3(\bfA - \bfB) +
          (\bfA-\bfB)^{T}\bfH_3\bfB \right)\right|\\
 & \leq & \onehalf\left(\sum_{i=1}^{3}|\bfA^{T}\bfH_{i}(\bfA - \bfB) +
        (\bfA-\bfB)^{T}\bfH_{i}\bfB|^2\right)^{1/2}\\
 & \leq & \onehalf\left(\sum_{i=1}^{3}|\bfA^{T}\bfH_{i}(\bfA - \bfB) +
        (\bfA-\bfB)^{T}\bfH_{i}\bfB|\right)\\
 &= & \onehalf\sum_{i=1}^{3}|\bfH_{i}|~\left(|\bfA| +
        |\bfB|\right)|\bfA - \bfB|\\
 & \leq & \tilde\beta \left(|\bfA| + |\bfB|\right)|\bfA - \bfB|.
\end{eqnarray*}
Since $|\bfH_{i}|= \sqrt4 =2$,
the proof  is completed  with
$\tilde\beta  =3$.

\begin{lemma}
Let $V(\bfx)$ be defined as in (\ref{Vol}).
Let $|\bfx_k - \bfx|_{1,4} \rightarrow 0$   where $\bfx,\bfx_k\in X$
and
$V(\bfx_k) = \omega_{0}$,
then
$$
\lim_{k\rightarrow \infty} V(\bfx_k) = V(\lim_{k\rightarrow \infty} \bfx_k )
= V(\bfx) = \omega_0.
$$
\label{XgVolClosed}
\end{lemma}
For each $k$, we have
\begin{eqnarray*}
|V(\bfx) - \omega_{0}| & = & |V(\bfx) - V(\bfx_k)|\\
        & = & \left|\onethird\int_{\Omega}(\bfx\cdot
        \adjugate_2\nabla\bfx - \bfx_k\cdot
        \adjugate_2\nabla\bfx_k) dA \right|.
\end{eqnarray*}
Adding
and subtracting $\bfx\cdot \adjugate_2\nabla\bfx_k$ to this last
expression and  applying Lemma~\ref{ADJlemma},
  we are led to
\begin{eqnarray}
|V(\bfx) - \omega_{0}| & = &
\left|\onethird\int_{\Omega}(\bfx-\bfx_k)\cdot
\adjugate_2\nabla\bfx_k +
\bfx\cdot (\adjugate_2\nabla\bfx - \adjugate_2\nabla\bfx_k) dA\right|\cr
        & \leq & \int_{\Omega}|\bfx-\bfx_k||\nabla\bfx_k|^{2} dA \ \ +\cr
    && \mbox{\hskip 1.5cm}
         \int_{\Omega} |\bfx|\left(|\nabla\bfx| +
        |\nabla\bfx_k|\right)|\nabla\bfx - \nabla\bfx_k| dA.
    \label{QQQ}
\end{eqnarray}
 Since $\bfx, \bfx-\bfx_k , \nabla\bfx - \nabla\bfx_k  \in L^{4}(\Omega,
\realR^3),$ $\nabla\bfx, \nabla\bfx_k \in L^{2}(\Omega,
\realR^3),$ and $|\nabla\bfx_k|^{2} \in L^{4/3}(\Omega,
\realR^3),$ we can apply the generalized H\"{o}lder inequality
to (\ref{QQQ})
and obtain
\begin{eqnarray}
|V(\bfx) - \omega_{0}| & \leq &
\left(\int_{\Omega}|\bfx-\bfx_k|^{4} dA
\right)^{1/4}\left(\int_{\Omega}|\nabla\bfx_k|^{\frac{8}{3}} dA
\right)^{3/4}  +   \cr
    & &  \left(\int_{\Omega} |\bfx|^4~dA\right)^{1/4}\left[
         \left(\int_{\Omega}|\nabla\bfx|^2~dA\right)^{1/2}
       +  \right. \cr
    & & \left.   \ \ \left(\int_{\Omega}|\nabla\bfx_k|^2
dA\right)^{1/2}\right]
\left(\int_{\Omega}|\nabla\bfx - \nabla\bfx_k|^4~dA\right)^{1/4}
\label{InEQlast}
\end{eqnarray}
Since  $\bfx_k \rightarrow \bfx \in W^{1,4}(\Omega, \realR^3)$ and
$\norm{\nabla\bfx_k}_{L^4}$ is uniformly
bounded,
from (\ref{InEQlast})
we see that
$|V(\bfx) - \omega_{0}| = 0$ and the proof
is complete.

\begin{remark}
{\rm
From Lemma~\ref{XgVolClosed}, we see that $ V(\bfx)-\omega_0 =0$
is a global constraint satisfying (\ref{globalCOND}) in
Lemma~\ref{XgClosed}.
When  $V(\bfx) -\omega_0=0$
is included
among  the  constraints,
$\bfg \le \bfzero$, we denote the corresponding closed subspace  of $X$ by
$X_\bfgo$ in order to emphasize the dependence on $\omega_0$.
}\end{remark}

\section{Existence Results}

\label{sec:4}
\setcounter{equation}{0}

In the following, we write $\bfu =(u,v)\in \Omega$.
The total energy  of the balloon system is taken to be
the sum of the hydrostatic pressure potential and the film strain energy,
\be
 \Etotal^*(\bfx, \nabla \bfx)  =  \int_\Omega \left(
W_f^*(\bfu, \bfx, \nabla\bfx)   +
               f_P(\bfx, \nabla \bfx)      \right)~dA.
 \label{totalE}
\ee
Combining (\ref{fPBounds}) and  (\ref{UpLowB}),
we find
\be
W_f -\onefourth\tau(1+\nu)
-  \left(\onehalf bR^2 + |p_0|R\right) |\bfF|^2
\le
  W_f^* +  f_P
\le
W_f  + \left( \onehalf bR^2 + |p_0|R\right) |\bfF|^2.
\label{UpLowBf}
\ee
Next, applying the appropriate upper and lower bounds on
$W_f$ established in (\ref{WfBounds}), we have
$$
\begin{array}{l}
\oneeighth \tau \left( \nu |\bfF|^4 - 2(1+\nu)|\bfF|^2 + 2(1 +  \nu)
 \right)
 -\onefourth\tau(1+\nu)-\left( \onehalf bR^2 + |p_0|R\right)   |\bfF|^2
\cr
\ \ \ \
\le   W_f^* +  f_P \le
 \oneeighth \tau |\bfF|^4 +  \onefourth \tau(1+\nu)
 + \left( \onehalf bR^2 + |p_0|R\right) |\bfF|^2.
 \end{array}
$$
Simplifying, we have
\be
\begin{array}{c}
\oneeighth \tau  \nu |\bfF|^4 -\onefourth \left(\tau(1+\nu) +
    2bR^2 +4|p_0|  R \right)|\bfF|^2 \cr
  \le W_f^* +  f_P\cr
\le
 \oneeighth \tau |\bfF|^4 +  \onefourth \tau(1+\nu)
 + (\onehalf bR^2 + |p_0|R)|\bfF|^2.
 \end{array}
 \label{BBB}
\ee Let $\kappa_1= \oneeighth \tau \nu$ and $\kappa_2 =
\onefourth(\tau(1+\nu) +2 bR^2 +4 |p_0|R)$ and apply
Lemma~\ref{CorINEQ}(i) to the first inequality in (\ref{BBB}).
Choose $\alpha_1$ so that $0< \alpha_1 < \oneeighth \tau \nu$ and
define $\nu_1$ by the relation $\alpha_1 = \oneeighth\tau\nu_1$. We
see that $0<\nu_1<\nu$ and by Lemma~\ref{CorINEQ}(i), we have
$$
\oneeighth \tau \nu_1 |\bfF|^4 - \rho_1
\le
\oneeighth \tau  \nu |\bfF|^4
  - \onefourth (\tau(1+\nu) + 2bR^2 + 4|p_0|R)|\bfF|^2,
$$
where
$$
{\displaystyle \rho_1 = \frac{ (\tau(1+\nu) + 2b R^2 + 4|p_0|R)^2}{8\tau(\nu-\nu_1)} }.
$$
Turning to the second inequality in  (\ref{BBB}),
and applying Part (ii) of Lemma~\ref{CorINEQ},
we choose
$\theta_1>1$, and
it follows that
$$
\oneeighth \tau |\bfF|^4
 + \left(\onehalf bR^2+|p_0|R\right) |\bfF|^2 \le
 \oneeighth\tau \theta_1 |\bfF|^4 +\rho_2,
$$
where
$$
\rho_2=  \frac{\left( bR^2+2|p_0|R\right)^2}{ 2\tau( \theta_1 -1)}.
$$
Thus,
$$
\oneeighth\tau \nu_1 |\bfF|^4 - \rho_1 \le
 W_f^* +  f_P  \le
\oneeighth\tau \theta_1 |\bfF|^4 + \rho_2 + \onefourth \tau(1+\nu),
$$
and
\be
\oneeighth\tau \nu_1 |\bfF|^4  \le
   W_f^* +  f_P  + \rho_1 \le
\oneeighth \tau \theta_1 |\bfF|^4 + \onefourth \tau(1+\nu) + \rho_1+\rho_2.
\label{ImportantIneq}
\ee
It follows that
\be
 f_T^*(\bfu,\bfx,\nabla\bfx ) = W_f^*(\bfu,\bfx,\nabla\bfx )
         + f_P(\bfu,\bfx,\nabla\bfx ) + \rho_1
\label{fTDEF}
\ee
is quasiconvex and satisfies some
additional properties
that we summarize in the following lemma.

\def\fstar{{f_T^*}}

\begin{lemma}
If $f_T^*(\bfu,\bfx,\nabla\bfx ) = W_f^*(\bfu,\bfx,\nabla\bfx )
         + f_P(\bfu,\bfx,\nabla\bfx ) + \rho_1$ and  $|\bfx| < R$, then
$\fstar$ is quasiconvex and
\begin{description}
\item[(i)]
${\displaystyle \gamma |\bfA|^4 \le   f_T^*(\bfu, \bfx,\bfA) \le
\delta( 1+  |\bfA|^4) }$
\item[(ii)]
$
 \left| \fstar(\bfu,\bfx, \bfA) - \fstar(\bfu,\bfy,
\bfB)\right|$ $ \le \beta \left(1 + |\bfx|^3 +  |\bfy|^3 +  |\bfA|^3 + \right.$
 \hfil\break \mbox{\space \hskip 6cm\space }  $\left.|\bfB|^3\right) \left( |\bfx - \bfy|
 +
 |\bfA - \bfB|
\right).$
\end{description}
\label{fstarProp}
\end{lemma}

Property {(i)} is (\ref{ImportantIneq}) with
$\delta = \max\lbrace \oneeighth\tau
\theta_1, \onefourth \tau(1+\nu) + \rho_1 + \rho_2 \rbrace$
and $\gamma = \oneeighth\tau \nu_1$.  We
will show Property~(ii) as follows.  Consider
\begin{eqnarray}
\fstar(\bfu,\bfx,\bfA) - \fstar(\bfu,\bfy,\bfB) & = &
 \left(W_f^*(\bfu,\bfx,\bfA ) - W_f^*(\bfu,\bfy,\bfB) \right)  \cr
&& \ \ \ + \ \ \left(f_P(\bfu,\bfx,\bfA ) - f_P(\bfu,\bfy,\bfB) \right).
\label{iiInEq}
\end{eqnarray}
Since
Lemma~\ref{appendixlemma} applies to $W_f^*$, there exists $\beta_1
> 0$ such that
\be
|W^*_f(\bfu,\bfx,\bfA) - W^*_f(\bfu,\bfy,\bfB)|  \le  \beta_1(1 +
|\bfA|^3 + |\bfB|^3)|\bfA - \bfB|.
\label{comb}
\ee
To complete the proof of Property~(ii),
$\pm[\onehalf
b(\bfy\cdot\bfk)^2 + p_0(\bfy\cdot\bfk)]\bfk\cdot\adjugate_2\bfA$
is added
to the last term   in~(\ref{iiInEq}). Using
 (\ref{fPdef})-(\ref{gPdef}),
we find \be f_P(\bfu,\bfx,\bfA ) - f_P(\bfu,\bfy,\bfB) = I_1 + I_2 +
I_3 + I_4, \ee where
\begin{eqnarray*}
I_1 &=& -\onehalf
b  \left[(\bfx\cdot\bfk)^2 -
 (\bfy\cdot\bfk)^2\right] \bfk\cdot\adjugate_2\bfA,  \label{I1} \cr
 I_2 &=&-p_0 \left[ (\bfx\cdot\bfk) -
(\bfy\cdot\bfk) \right] \bfk\cdot\adjugate_2\bfA,  \label{I2} \cr
I_3 &=&-\onehalf
b  \left[ (\bfy\cdot\bfk)^2  \right] \bfk \cdot \left( \adjugate_2\bfA -
\adjugate_2\bfB\right),  \label{I3}\cr
I_4 &=&  -p_0  \left[ \bfy\cdot\bfk \right] \bfk \cdot \left(\adjugate_2\bfA -
\adjugate_2\bfB\right). \label{I4}
\end{eqnarray*}
By   Young's inequality,
we have $ 3|\bfx|~|\bfA|^2 \le  |\bfx|^3 + 2|\bfA|^3 $,
 and $ 3|\bfy|~|\bfA|^2 \le  |\bfy|^3 + 2|\bfA|^3$.
 Using these estimates
and
Lemma~\ref{ADJlemma} with $\tilde\beta=3$, we find
 \begin{eqnarray}
|I_1| & \le & \onehalf b |\bfx - \bfy|
\left(  |\bfx|+|\bfy| \right) | \adjugate_2 \bfA|\cr
& \le& \threehalves b |\bfx-\bfy|~\left( |\bfx|+ |\bfy|\right) |\bfA|^2\cr
&\le &  \onehalf b \left( |\bfx|^3 + |\bfy|^3 + 4|\bfA|^3 \right) |\bfx - \bfy|\cr
&\le &  2 b \left(
1+ |\bfx|^3 + |\bfy|^3 +  |\bfA|^3 + |\bfB|^3\right) |\bfx - \bfy|.
\label{I1f}
 \end{eqnarray}
 Applying
Lemma~\ref{ADJlemma} with $\tilde\beta=3$
and Young's inequality with  $3|\bfA|^2 \le 1 + 2 |\bfA|^3$, we find
 \begin{eqnarray}
|I_2| & \le &  |p_0|~|\bfx - \bfy|~|\adjugate_2 \bfA| \cr
      & \le &  3|p_0|~|\bfx - \bfy|~|\bfA|^2 \cr
      & \le &  |p_0|~\left(  1 + 2 |\bfA|^3\right)~|\bfx - \bfy| \cr
      & \le &  2|p_0|~\left( 1+ |\bfx|^3 + |\bfy|^3 +  |\bfA|^3
      + |\bfB|^3 \right)~|\bfx - \bfy|.
\label{I2f}
 \end{eqnarray}
Applying Lemma~\ref{ADJlemma} with $\tilde\beta=3$
and Young's inequality with  $3|\bfy|^2 |\bfA| \le  2 |\bfy|^3  +  |\bfA|^3$
and
$3|\bfy|^2 |\bfB| \le  2 |\bfy|^3  +  |\bfB|^3$,
 we find
 \begin{eqnarray}
|I_3| & \le & \onehalf b |\bfy|^2 |\adjugate_2 \bfA - \adjugate_2 \bfB| \cr
      & \le & \threehalves   b |\bfy|^2 \left( |\bfA| + |\bfB| \right) | \bfA - \bfB| \cr
      &\le & \onehalf \left( 4 |\bfy|^3 + |\bfA|^3 + |\bfB|^3  \right)~|\bfA - \bfB| \cr
      &\le & 2 b \left(   1+ |\bfx|^3 + |\bfy|^3 +  |\bfA|^3 + |\bfB|^3
             \right)~|\bfA - \bfB|.
\label{I3f}
 \end{eqnarray}
 Applying Lemma~\ref{ADJlemma} with $\tilde\beta=3$
and  the estimates
$|\bfy|~|\bfA| \le \onethird + \twothirds |\bfy|^{3/2} |\bfA|^{3/2}
\le \onethird( 1 + |\bfy|^3 + |\bfA|^3 )$
leads  to
 \begin{eqnarray}
|I_4| & \le & |p_0|~|\bfy|~|\adjugate_2 \bfA - \adjugate_2 \bfB|\cr
      & \le &  3|p_0|~|\bfy|~\left( |\bfA| + |\bfB| \right)  |\bfA - \bfB| \cr
      & \le &  |p_0|~\left(
      2+ 2|\bfy|^3 + |\bfA|^3 + |\bfB|^3 \right)~|\bfA - \bfB| \cr
      & \le & 2|p_0|~\left(  1+ |\bfx|^3 + |\bfy|^3 +  |\bfA|^3 + |\bfB|^3
                  \right)~|\bfA - \bfB|.
\label{I4f}
 \end{eqnarray}
Combining (\ref{I1f})-(\ref{I4f}), we find
\be
\begin{array}{ccc}
|f_P(\bfu,\bfx,\bfA) - f_P(\bfu,\bfy,\bfB)| \cr
\multicolumn{3}{l}{\ \ \mbox{\hskip 0.5cm} \le 2(b+|p_0|)
\left(1 + |\bfx|^3 +  |\bfy|^3 +  |\bfA|^3 +
|\bfB|^3\right)\left(|\bfx - \bfy| + |\bfA - \bfB|\right).}\\
\end{array}
\label{almostLast}
\ee
Combining (\ref{iiInEq}),
(\ref{comb}),
and (\ref{almostLast}), we see
that Property~(ii) follows with
 with $\beta = 2\max\{\beta_1, 2(b+|p_0|)\}$.

Since all of the forces are conservative, equilibrium is achieved
at a minimum of the energy functional
over the class of functions  $X_\bfg$.
The constant $\rho_1$ serves only to increase
$W_f^*+f_P$ by a constant, so the minimizer of
$W_f^*+f_P$ is the same as the minimizer of $W_f^*+f_P+\rho_1$.
We state  two
existence
results for the equilibrium
shape of an inflated balloon.

\begin{theorem}{Closed balloon system.}
If  $|\bfx|\le R$, then
\begin{eqnarray*}
(P_\bfgo)\mbox{\hskip 1cm}
&
{\displaystyle
\inf \left\lbrace
I(\bfx) = \int_\Omega f_T^*(\bfu, \bfx,\nabla\bfx)~dA\ | \ \bfx \in  X_\bfgo\
\right\rbrace\mbox{\hskip 4cm}
}
\end{eqnarray*}
admits at least one solution.
\label{openTheorem}
\end{theorem}
Since we can always find an $\bfx \in \bfx_0 + W_0^{1,4}(\Omega;
\realR^3)$ parameterizing a spherical cap, we know that
$\calD_\bfgo$ is nonempty. In order to complete the proof,
  we need to show
 that $\fstar$ is quasiconvex and
 verify  H1~(i)-(iii)    in Theorem~\ref{appendixtheorem}.
Lemma~\ref{fstarProp} shows that $\fstar$ is quasiconvex and
satisfies
H1~(i)-(ii) if $|\bfx| \le R$.
H1~(iii) is
satisfied because $\fstar$ is independent of $\bfu$ and so $\eta\equiv 0$.
Lemma~\ref{XgVolClosed} establishes
that $X_\bfgo$ is closed.
It follows that
$(P_\bfgo)$ has at least one solution.

\begin{remark}
{\rm
Eq.~(\ref{Archimedes}) guarantees that
 the values of $b$  and $\omega_0$ are  sufficient
to lift the balloon system. If $b=b_d$ and $\omega_0=\omega_{0,d}$
then  the balloon envelope is fully deployed,
and we  are justified in ignoring
the
contribution of the film weight in the stress
analysis.}
\end{remark}
\begin{remark}
{\rm
In general, one should not expect a unique  solution
of $(P_\bfgo)$, especially if the balloon is not fully inflated.
For example, if the lift generated by the gas bubble
is not  large enough to lift the balloon or the gas is  compressed sufficiently
(i.e., $\omega_0 < \omega_{0,d}$),
then one  expects to find  multiple local equilibria
with slack regions where a unique
equilibrium
configuration may not
exist.
For an ascending large scientific
balloon,
it is not uncommon
to observe balloon shapes with a nearly periodic
lobe pattern surrounding the
gas bubble. These shapes are characterized by significant regions of
excess folded material hanging beneath the gas bubble.
Folded material can be handled by introducing constraints
in
the form covered by Lemma~\ref{XgClosed}
and
the theory developed in this paper
applies to shapes with  a period lobe pattern.}
\end{remark}

\begin{theorem}{Open balloon system.}
If $|\bfx| \le R$, then
\begin{eqnarray*}
(P_\bfg)\mbox{\hskip 2cm}
&
{\displaystyle
\inf \left\lbrace
I(\bfx) = \int_\Omega f_T^*(\bfu, \bfx,\nabla\bfx)~dA\ | \ \bfx \in  X_\bfg
\right\rbrace
}
\mbox{\hskip 4cm}
\end{eqnarray*}
 admits
 at least one solution.
\label{closedTheorem}
\end{theorem}

Since the subspace $X_\bfg$ is closed  by Lemma~\ref{XgClosed},
the  proof of
Theorem~\ref{closedTheorem} is the same as
the proof of Theorem~\ref{openTheorem} with the space
$X_\bfgo$ replaced by $X_\bfg$.

\section{Numerical results on strained balloon shapes}

\label{sec:5}
 \setcounter{equation}{0}

In this section, we present numerical solutions
of our model.
We consider a zero-pressure natural shape  design and
a pumpkin design. We assume the balloon film is
$32\mu$m polyethylene and the balloons have no
caps. We use the same tendon weight density
and
a suspended payload of 4000~N
for both
balloons  in shape finding.
The specific buoyancy at float is $b=0.068$~N/m$^3$
and corresponds to an altitude of roughly  35~km.
For the pumpkin balloon,
$r_B=0.785$~m and $p_0=200$~Pa.
See
Table~\ref{ShapeFindingParm} for a summary of
shape finding parameters that were used.
We used mechanical properties
that were determined by
the  Balloon Lab at NASA's Wallops Flight Facility
for
$0.8$~mil polyethylene film
at room temperature (this  is comparable to the temperature
at float
 during nominal daylight conditions).
For an isotropic three-dimensional
material, one can show that Poisson's ratio
is
between $0$ and $0.5$ (see, \cite[p. 129]{ciarletI}). However, in lab experiments
using  cylinder tests, one finds that for a thin
polyethylene
film,
Poisson's ratio is greater than  $0.5$ (see, e.g., \cite{blandino}).
We extrapolated the findings for 0.8 mil to 1.5 mil film
to estimate  Young's modulus for our simulations.
Mechanical properties are summarized
in
Table~\ref{MechProp}.
\begin{table}
\begin{center}
\footnotesize
\begin{tabular}{lcc}
\hline
Description &Variable & Value\cr
\hline
Number of gores         &  $n_g$ & 200\cr
Buoyancy (N/m$^3$)           &$ b$  & 0.068\cr
Tendon weight density (N/m)&$w_t$ & 0.094 \cr
Film weight  density (N/m$^2$)&$w_f$ & 0.344 \cr
Payload  (N)     & $ L$ &  4000\cr
End-plate diameter (m)       & $ d_1$  & 1.32   \cr
\hline
\end{tabular}
\end{center}
\caption{\captionsize Shape finding  parameters. For pumpkin balloons
we assumed $r_B= 0.785$~m  and $p_0=200$~Pa.
For ZPNS balloons, $p_0=0$~Pa.
Tendon weight and film weight are included in the
shape finding process.}
\label{ShapeFindingParm}
\end{table}

\begin{table}
\footnotesize
\begin{center}\captionsize
\begin{tabular}{lcc}
\hline
Description & Variable & Value\cr
\hline
Film Young's modulus (MPa) & $E_f$ &   404.2 \cr
Film Poisson ratio    &   $\nu_f$& 0.825\cr
Film thickness   ($\mu$m) & $t$      & 32 \cr
\hline
\end{tabular}
\end{center}
\caption{\captionsize  Mechanical properties of  polyethylene film.}
\label{MechProp}
\end{table}

Next, we summarize our findings when we numerically solved
Problems~$(P_\bfg)$  and $(P_\bfgo)$.
In all cases, we assumed a cyclically symmetric shape and so
we needed to solve for one-half a gore.
We divided a half-gore into three strips
with 100 triangles per
 strip.
We used our own finite element code written in {\sf Matlab} to
compute the energies and constraints as described in Sections~3-5.
We then used \mbox{\sf fmincon} from {\sf Matlab's Optimization
ToolBox} to solve a constrained minimization problem where $\fstar$
was the objective function and $\bfg \le 0$ were the constraints.
Gradients for $\fstar, \bfg$, and the Hessian of $\fstar$ were
computed analytically. When  a nonlinear constraint (i.e., volume or
tendon constraint) is  imposed, we used a ``medium-scale''
sequential quadratic programming quasi-Newton's method; when there
were only linear constraints, a ``large-scale'' projected
trust-region Newton's method was used (see, \cite{matlab} for
further details on these algorithms).

We discretized the three-dimensional surface produced by
the shape finding process and used this as an initial
guess to start the solution process.
For each of the scenarios  considered, we present
plots of the ``averaged'' principal strains and
``averaged'' principal stress resultants
for each adjacent pair of rectangles in a strip.
If a meridional strip has $2M$ triangles
(numbered from top to bottom) and
$\mu_{ 2i-1,1}$ and $\mu_{2i,1}$ are the meridional  principle
stress resultants, then we plot
$\onehalf ( \mu_{ 2i-1,1} +   \mu_{2i,1})$ for $i=1, \dots M$.
Note, the directions for $\mu_{ 2i-1,1}$ and $\mu_{2i,1}$ need not
coincide, but this measure seems to work satisfactorily for
summarizing the data.
A similar convention is followed for averaging
the circumferential stress resultants and
the principal strains.

\subsection{ Zero-pressure natural shape balloon.}

\label{case:1}

We considered two scenarios involving a ZPNS balloon.

\begin{description}
\item{(a)} {\em Closed system $V(\bfx) = \omega_0$.}
In this case,
we considered a closed system with
a target volume $\omega_0=137,023$~m$^3$.
The averaged principal strains are plotted in
Figure~\ref{volZPNS}(a) and the averaged
principal stress resultants
are
plotted in Figure~\ref{volZPNS}(b).
The $\delta_1$ and $\delta_2$ directions
corresponds roughly to
meridional and
circumferential, respectively.
We found that the maximum principal strains were less
than $0.03\%$. The meridional stress resultants are
plotted in the top graph in Figure~\ref{volZPNS}(b)
and increase monotonically
from the bottom of the gore to its top.
From
 Figure~\ref{volZPNS}(b), we see that
near the top one-third, the film is in a biaxial state
and in the middle one-third it is in a uniaxial or near uniaxial
state  (i.e., $\mu_{2,i} \approx 0$).

\item{(b)} {\em Open system $p_0=0$.}
In Section~\ref{case:1}(b), we considered the same conditions
as those
in Section~\ref{case:1}(a), except we dropped the volume constraint
and set $p_0=0$. The results were nearly identical to those
obtained in Section~\ref{case:1}(a) and for this reason plots
were not included.
We found that the difference between the volumes in
Section~\ref{case:1}(a) and (b) was  less than $10^{-5}$~m$^3$.

\end{description}

In the following, principal strains and stress resultants
are to be interpreted as
`averaged' principal strains and
`averaged' principal
stress resultants, respectively.

\subsection{ Pumpkin balloon without tendons.}
\label{case:2}

To get a sense of the important
role that tendons
play in a pumpkin balloon,
we computed a solution for a pumpkin balloon
without tendons.
Since it was clear that a single layer of 32~$\mu$m film
was not sufficiently strong, we quadrupled the  thickness for
this numerical experiment.
From
Figure~\ref{NoTendonsNoVol}(a) we see that
the strains are nearly 7.0\%
and the principal stress resultants are quite large
(a maximum of  8~kN/m).
It is interesting to note that near the  ends,
the deformed gore is in a state of biaxial tension,
where in the middle of the gore the tension is
uniaxial and parallel to the gore length.
In reality,
the film is a  viscoelastic
material and under the right conditions  may be able to
maintain the integrity of the
gas barrier at
such a
high strain.  However,
the peaks in  Figure~\ref{NoTendonsNoVol}
suggest
this loading scenario should be avoided.

\subsection{Pumpkin balloon with tendons.}
\label{case:3}

In this case, we considered a
single layer of 32~$\mu$m polyethylene
film and added inextensible  tendons
via a constraint (i.e., (\ref{tCON})).
In Figures~\ref{withTendons}(a)-(b), we present
the principal strains
and principal stress resultants.
We see that when tendons are added,
the maximum principal strains are reduced from
  7.0\% in Section~\ref{case:2}
 to 1.5\%.
The  maximum principal
stress resultants
are reduced from
 8~kN/m in Section~\ref{case:2} to 225~N/m.
It is interesting to note from the plot of the
circumferential stress resultants (bottom graph of
Figure~\ref{withTendons}(b)), the film is in a state of biaxial
stress along the center of the gore. But near the equator
and the sub-tropic regions, as one gets
closer to the tendon,
the film is in a uniaxial state. This makes sense because
one would expect there to be wrinkling in these
locations.

\subsection{Pumpkin balloon with shortened tendons}
\label{case:4}
In Figures~\ref{withTendonsShort}(a)-(b), we plot
the principal strains
and principal stress resultants when shortened tendons are utilized.
Basically, we reduced
the tendon length in Section~\ref{case:3} by~2\%. Overall,
we see that the principal  strains and principal stress resultants
are reduced. However,
the maximum  principal  strains and principal stress resultants
are not reduced by a significant amount.
Comparing Figure~\ref{withTendons}(b) and
Figure~\ref{withTendonsShort}(b), we see
more
wrinkling in Figure~\ref{withTendonsShort}(b), especially in
the area adjacent
to the tendons.

While there are many other factors that go into the design of a gore,
a   balloon of the type discussed
in Section~\ref{case:4}
is not  very efficient.  Ideally, one would like to have the gore
uniformly loaded. The presence of wrinkles in the fully inflated shape
suggests that there is excess material that is doing little to carry
its share of the load. Furthermore, it has been demonstrated
that
too much tendon shortening can be detrimental,  leading
to
instability  in the strained cyclically symmetric
float shape  and
impeding proper deployment (see, \cite{JAaustin}-\cite{ASRparis}).

%

\section{Conclusions}
\label{sec:6}
 \setcounter{equation}{0}

Motivated by the problem of a large scientific balloon at
float altitude, we
presented a mathematical model for
a strained inflated  wrinkled membrane  loaded by hydrostatic pressure and
constrained by
tendons.
We assumed the balloon is constructed of a   thin
linearly elastic isotropic
material.
Load tendons are included through the use of constraints.
Using direct methods in the calculus of variations
we established  rigorous existence theorems under
 general conditions.
We computed numerical solutions based on our model and
estimated the principal
strains and principal stress resultants under
nominal loading conditions at float altitude.
Our theoretical results establish a solid foundation
for our mathematical model and  affirms
the use of  numerical
computations to estimate film stress resultants, information
that is valuable to the balloon designer.
An efficient balloon design keeps weight
to a minimum,
but the balloon  must be of sufficient strength to operate
safely over its service life.
Analytical and computational
tools are demonstrated here that can
help
the balloon designer balance the competing factors of
a stronger, albeit heavier balloon versus
a lighter more efficient~one.

\centerline{\bf Acknowledgement}
The first two authors  were  supported   by NASA Grant NAG5-5353.

\vfil\eject

\appendix

\newtheorem{lemmaAPP}{Lemma}[section]
\newtheorem{theoremAPP}{Theorem}[section]

\section{Direct methods in the calculus of variations}

\setcounter{equation}{0}
\label{sec:A}

For the convenience of the reader,
we record two results
by Dacorogna.
\begin{lemmaAPP}{ Dacorogna \cite[Lemma 2.2, p.~156]{dacorogna} }
Let $f: \realR^n \rightarrow \realR$ be convex in each variable
and let
$$
|f(x)| \le \alpha\left( 1 + |x|^p\right)
$$
for every $x \in \realR^n$ and where $\alpha \ge 0, p \ge 1$. Then there exists
$\beta\ge 0$ such that
$$
|f(x) - f(y)| \le \beta \left( 1 + |x|^{p-1} + |y|^{p-1} \right)|x-y|
$$
for every $x, y\in \realR^n$.
\label{appendixlemma}
\end{lemmaAPP}
\begin{theoremAPP}{Dacorogna \cite[Theorem 2.9, p. 180]{dacorogna}}
Let $\Omega \subset \realR^n$ be a bounded open set. Let
$f:\Omega\times \realR^m\times \realR^{nm} \rightarrow \realR$ be
continuous and quasiconvex and satisfying
\begin{description}
\item[(i)]
${\displaystyle
\gamma |\bfA|^p \le  f(\bfu, \bfx,\bfA)
\le  \delta( 1+ |\bfx|^p + |\bfA|^p)
}$;
\item[(ii)]
$ \left|
f(\bfu,\bfx, \bfA) - f(\bfu,\bfy, \bfB)\right|
\le
\beta
\left(1 + |\bfx|^{p-1} + |\bfy|^{p-1}   \right. \ +$ \hfil\break
\mbox{\hskip 4cm}$ \left. \ \ |\bfA|^{p-1} +  |\bfB|^{p-1}\right)$
$\left(
|\bfx - \bfy| + |\bfA - \bfB|
\right);
$
\item[(iii)]
${\displaystyle
\left|
f(\bfu,\bfx, \bfA) - f(\bfv,\bfx, \bfA)\right|
\le
\eta(|\bfu-\bfv|) \left(
1 + |\bfx|^p + |\bfA|^p
\right)
}$
where
$\eta$ is a continuous increasing  function with $\eta(0)=0$
and $ p>1, \alpha, \beta,\gamma>0$ are  constants.
\end{description}
Let
\be
(P)\mbox{\hskip 0.3cm}
\inf \left\lbrace
I(\bfx) = \int_\Omega f( \bfu, \bfx(\bfu), \nabla\bfx(\bfu))~dA \ : \
\bfx \in \bfx_0 + W^{1,p}_0(\Omega; \realR^m)
\right\rbrace
\label{APP}
\ee
then $(P)$ admits at least one solution.
\label{appendixtheorem}
\end{theoremAPP}

\twocolumn

\begin{figure}
\begin{center}\captionsize
\begin{tabular}{c}
\mbox{\psfig{figure=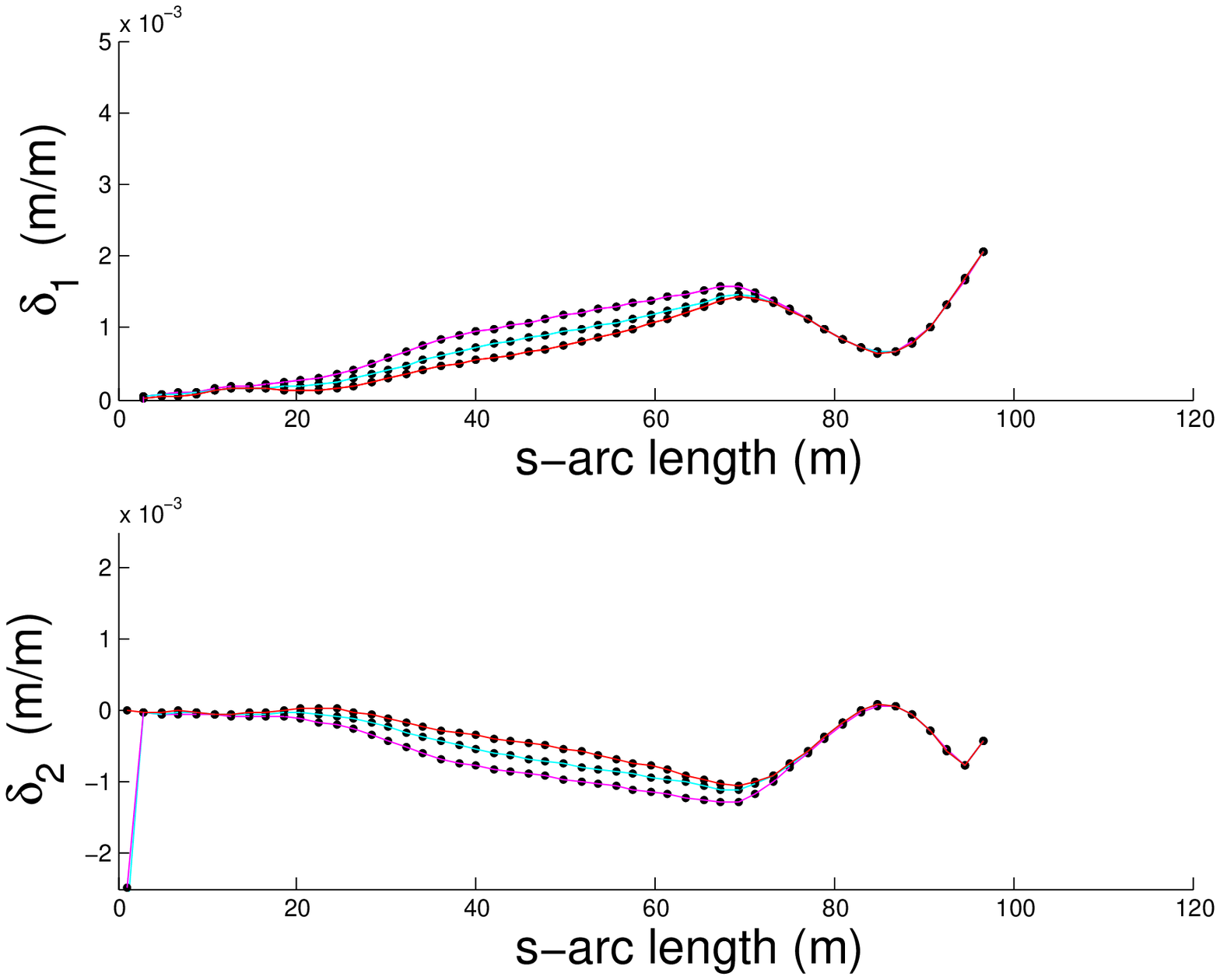,height=\FSIZEs }}
\cr
{\captionsize (a)~Averaged principal strains}  \cr
{\mbox{\psfig{figure=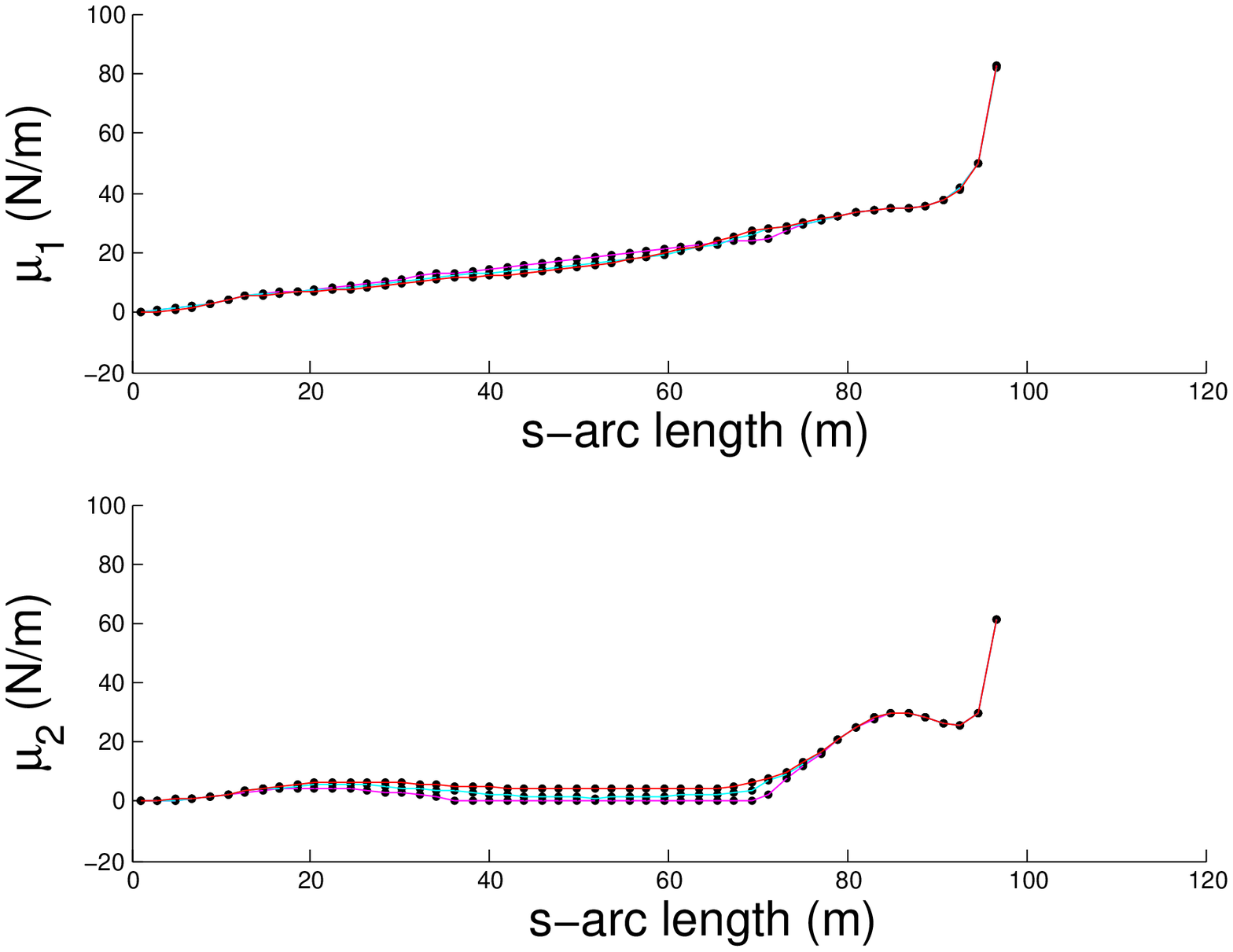,height=\FSIZEs }}}\cr
{\captionsize (b)~Averaged principal stress resultants}
\end{tabular}
\end{center}
\caption{\captionsize Zero-pressure natural shape balloon with
tendon and volume constraints.
}
\label{volZPNS}
\end{figure}

\begin{figure}
\begin{center}\captionsize
\begin{tabular}{c}
\mbox{\psfig{figure=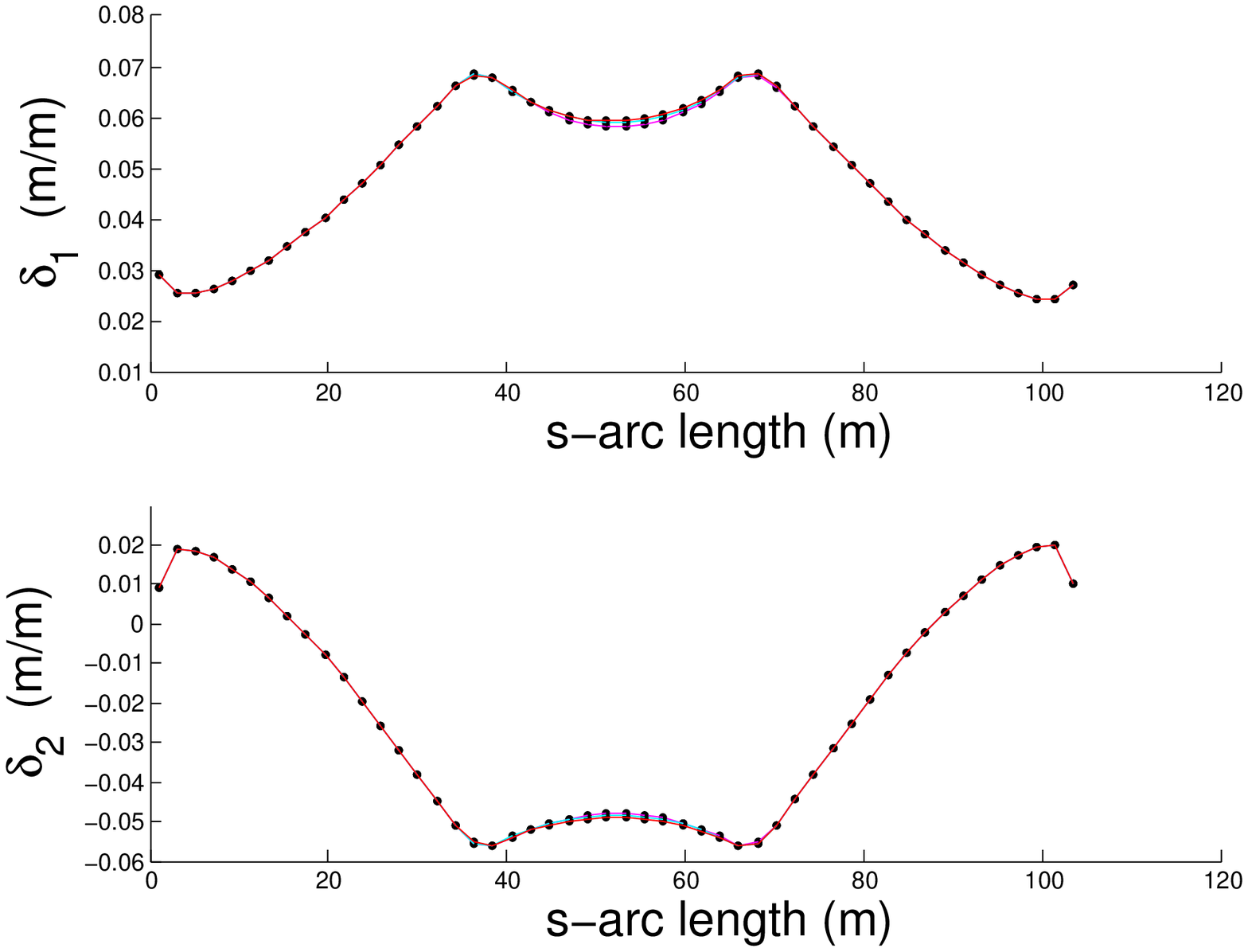,height=\FSIZEs }}
\cr
{\captionsize (a)~Averaged principal strains}  \cr
 {\mbox{\psfig{figure=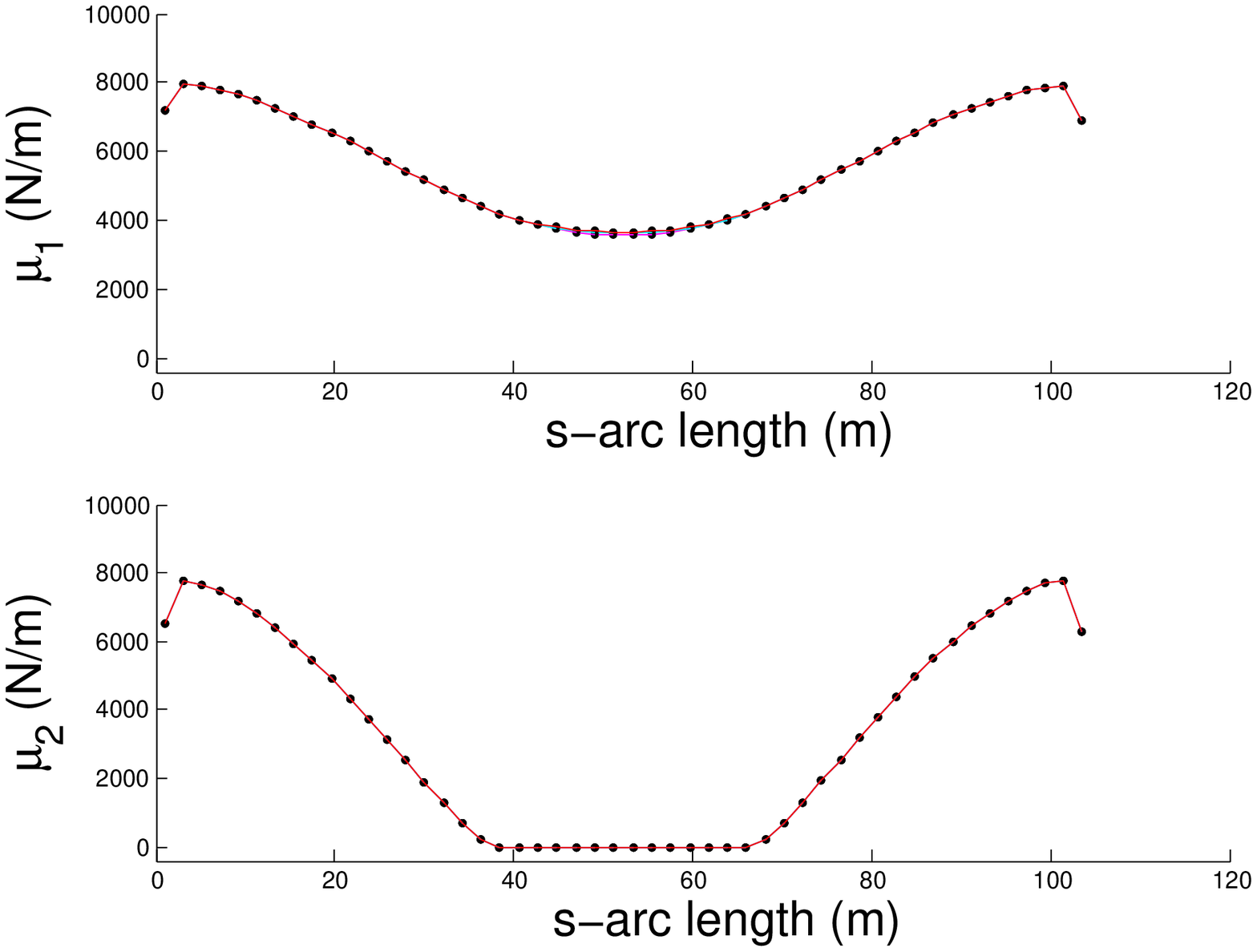,height=\FSIZEs}}}\cr
{\captionsize (b)~Averaged principal stress resultants}
\end{tabular}
\end{center}
\caption{\captionsize Pumpkin balloon without
tendon constraints and nominal thickness quadrupled.
}
\label{NoTendonsNoVol}
\end{figure}

\begin{figure}
\begin{center}\captionsize
\begin{tabular}{c}
\mbox{\psfig{figure=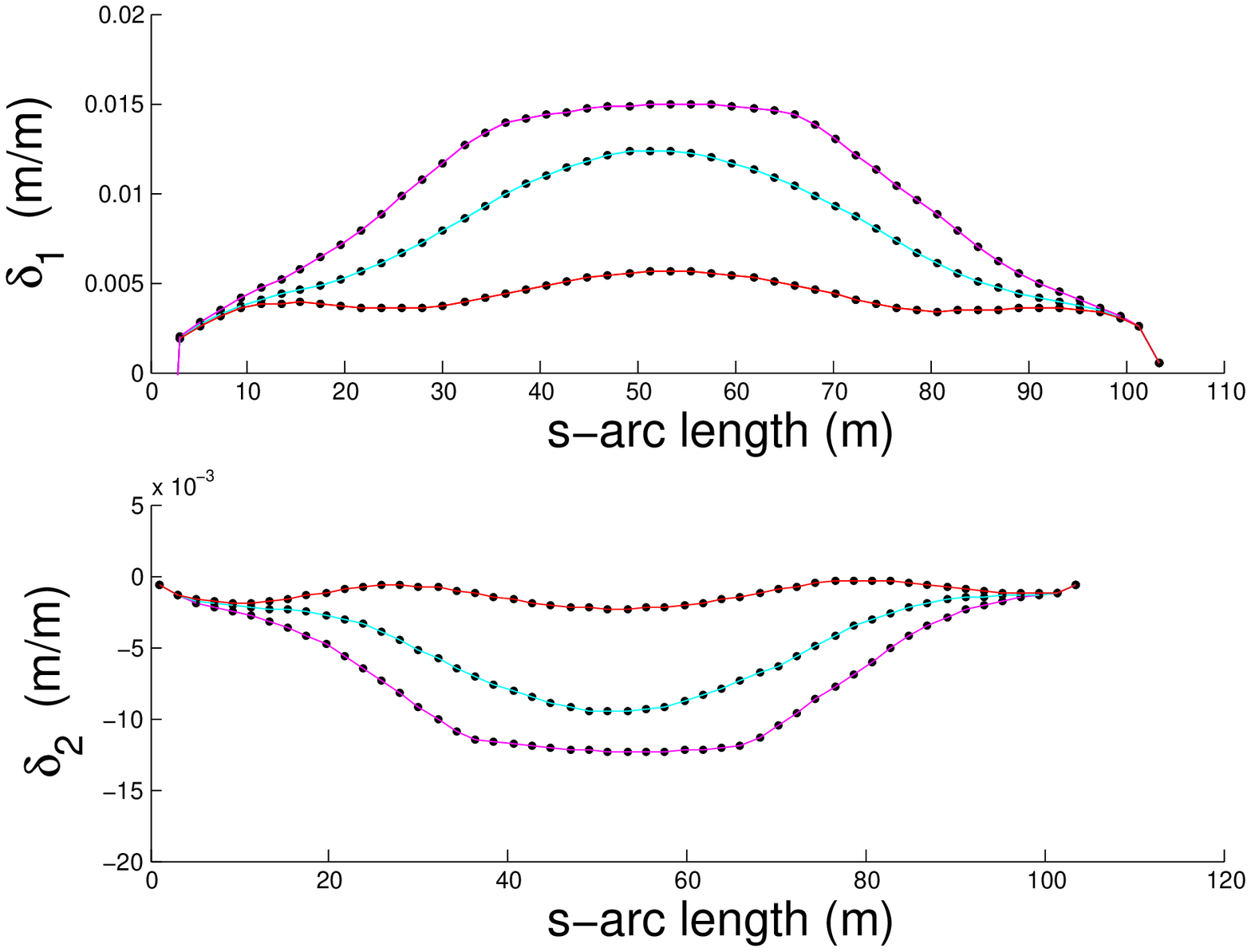,height=\FSIZEs }}
\cr
{\captionsize (b)~Averaged principal stress resultants}\cr
{\mbox{\psfig{figure=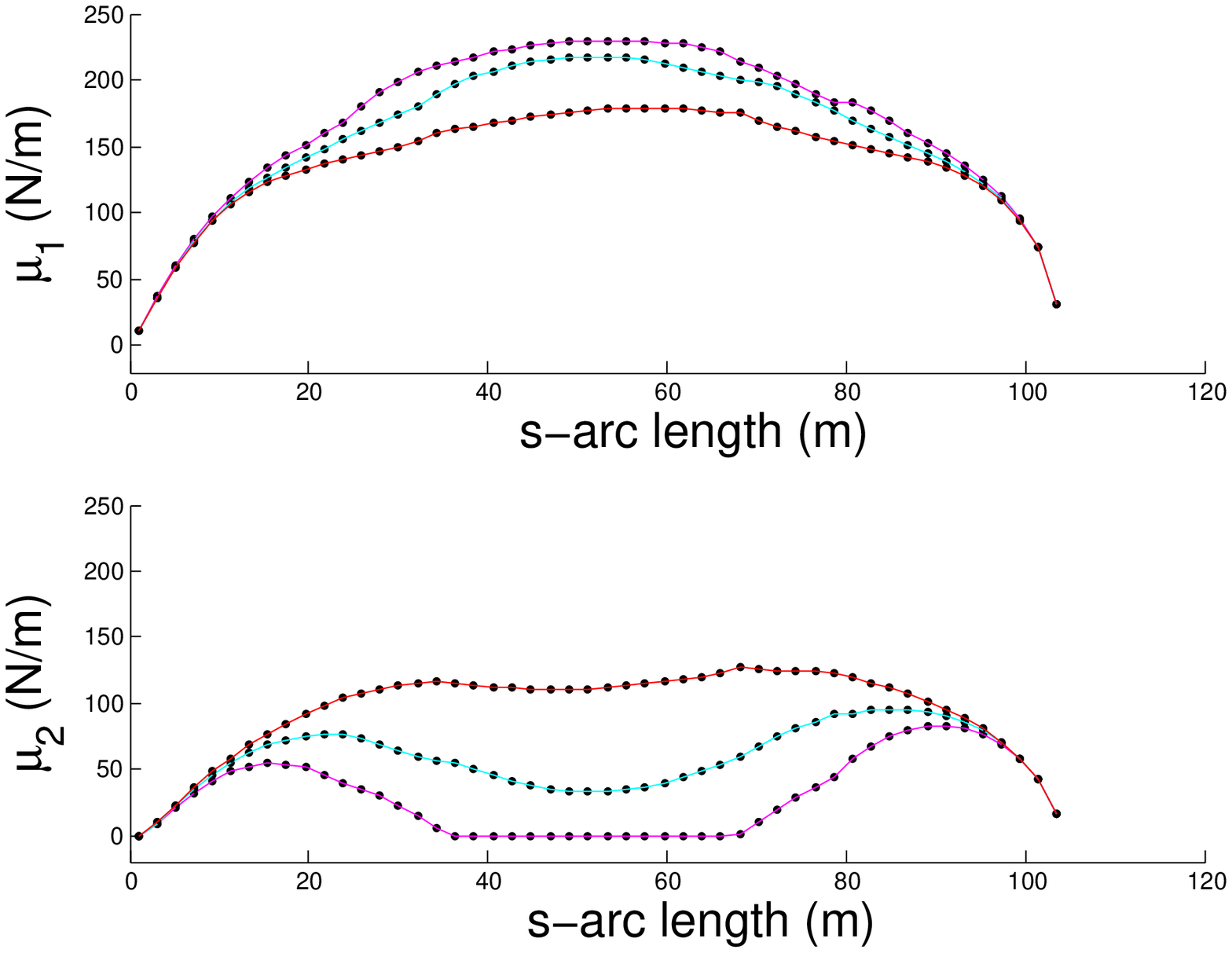,height=\FSIZEs }}}\cr
{\captionsize (b)~Averaged principal stress resultants}\cr
\end{tabular}
\end{center}
\caption{\captionsize Pumpkin balloon with
tendon constraints. }
\label{withTendons}
\end{figure}

\begin{figure}
\begin{center}\captionsize
\begin{tabular}{c}
\mbox{\psfig{figure=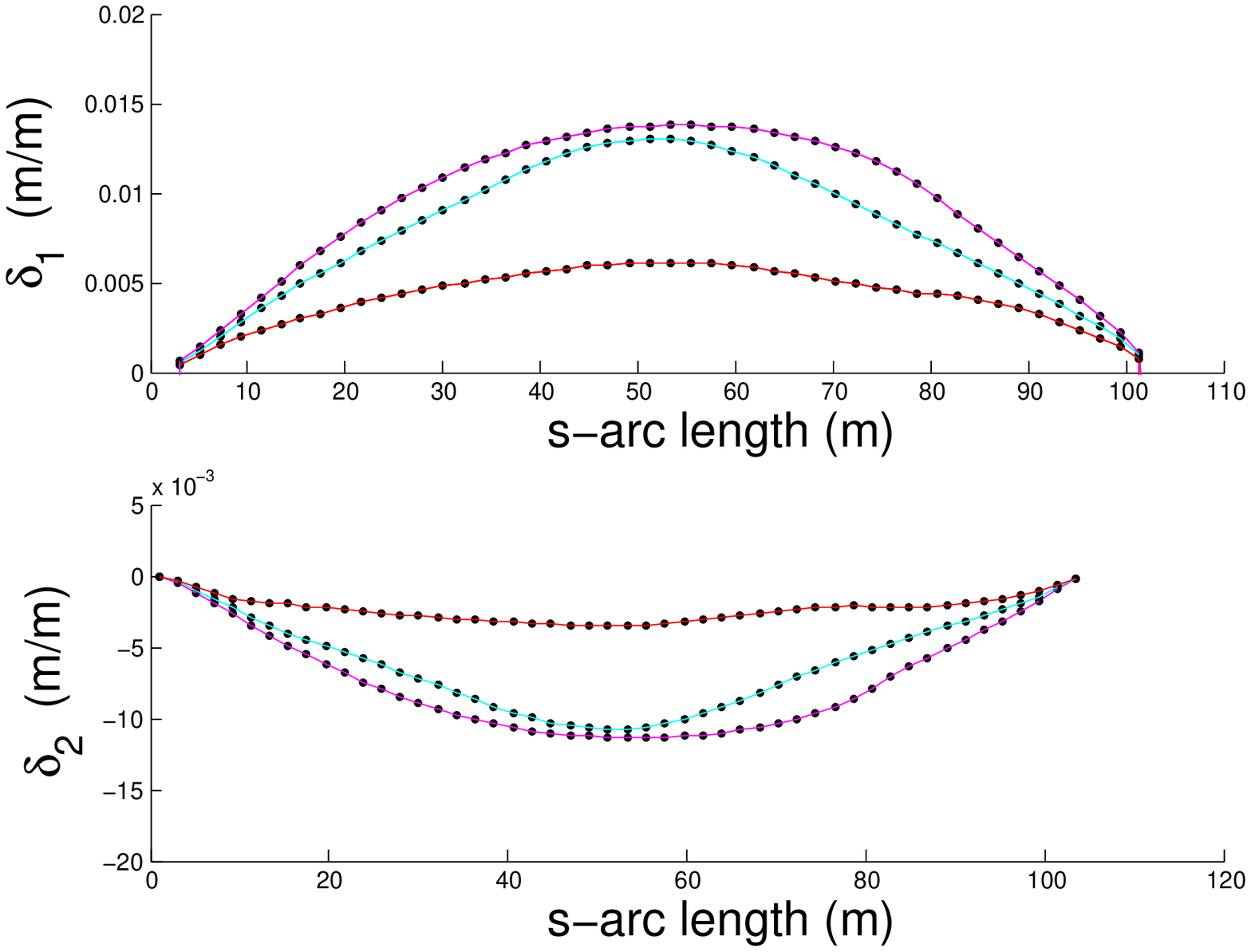,height=\FSIZEs }}
\cr
{\captionsize (a)~Averaged principal strains}  \cr
{\mbox{\psfig{figure=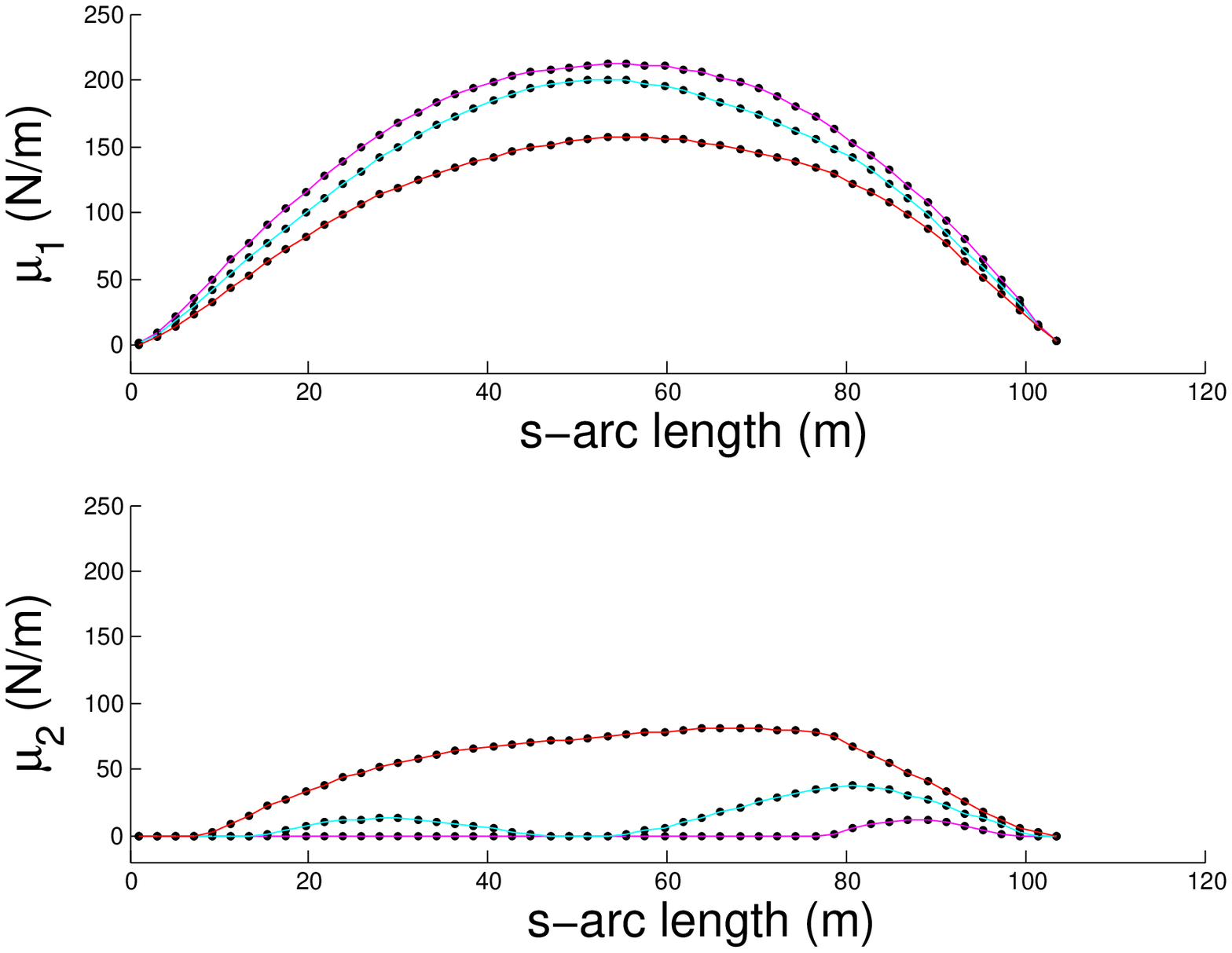,height=\FSIZEs }}}\cr
{\captionsize (b)~Averaged principal stress resultants}\cr
\end{tabular}
\end{center}
\caption{\captionsize Pumpkin
balloon with tendons shortened by
2\%.}
\label{withTendonsShort}
\end{figure}
\onecolumn



\begin{thebibliography}{}

\bibitem{umn} \textsc{ Anon.}:
\textit{Research Development in the Field of High Altitude Plastic Balloons},
NONR--710(01a) Reports,
Department of Physics, University of Minnesota,
Minneapolis, Minnesota,
1951--1956.

\bibitem{Ant} \textsc{S. S.  Antman}: \textit{Nonlinear Problems of Elasticity},
Springer-Verlag, New York, 1995.

\bibitem{aubin}
\textsc{T. Aubin}: \textit{Nonlinear analysis on Manifolds. Monge-Amp\'ere
Equations}, Springer-Verlag, New York, 1982.


\bibitem{BaginskiSIAP} \textsc{F. E. Baginski}:
\textit{On the design and analysis of inflated membranes: natural shape
and pumpkin shaped balloons},
SIAM J. Appl. Math. \textbf{65} No.~3  838-857 (2005).


\bibitem{JAaustin} \textsc{F.  Baginski, K. Brakke \& W. W. Schur}:
\textit{
Stability of cyclically
symmetric strained pumpkin balloon configurations and the formation of
      undesired equilibria},
      accepted for publication in the AIAA Journal of
      Aircraft.

\bibitem{JAcrystalcity}
\textsc{F. Baginski, K. Brakke \& W. W. Schur}:
\textit{Unstable cyclically symmetric and
stable asymmetric  configurations of a pumpkin balloon,}
 accepted for publication in the
AIAA Journal of Aircraft.

\bibitem{ASRparis} \textsc{F.  Baginski, K. Brakke \&  W. W. Schur}:
\textit{Cleft formation in pumpkin balloons},
Adv. Space Res.,
{\bf 37} 2070-2081 (2006).

\bibitem{PalmSprings} \textsc{F.  Baginski \&  W. W. Schur}:
\textit{Undesired equilibria of self-deploying pneumatic envelopes},
AIAA Journal of Aircraft,  \textbf{42}  No. 6 1639-1642  (2005).

\bibitem{BagSch} \textsc{F. Baginski \& W. Schur}:
\textit{
Structural analysis of pneumatic envelopes: A variational formulation
and optimization-based solution process},
AIAA Journal, \textbf{41} No. 2  304-311 (2003).



\bibitem{ball}  \textsc{J. M. Ball}:
\textit{\em Convexity conditions and existence theorems in
nonlinear elasticity},
Arch. Ration. Mech. Anal. \textbf{63} 337-403 (1977).


\bibitem{blandino} \textsc{J. Blandino,
J. Sterling, F. Baginski, E. Steadman, J. T. Black  \& R. Pappa}:
\textit{ Optical strain measurement of an inflated
cylinder using photogrammetry with application to
scientific balloons},
AIAA-2004-1500,
45th AIAA/ASME/ASCE/AHS/ASC Structures, Structural Dynamics,
and Materials Conference,
19-22 April 2004, Palm Springs, CA.


\bibitem{ciarletIII} \textsc{P. G. Ciarlet}:
\textit{Mathematical Elasticity Volume III:
Theory of Shells}, North-Holland, New York, 2000.

\bibitem{ciarletI} \textsc{P. G. Ciarlet}:
\textit{Mathematical Elasticity Volume I:
Three-Dimensional Elasticity}, North-Holland, New York, 1988.


\bibitem{matlab} \textsc{T. Coleman, M. A. Branch \& A. Grace}:
\textit{Optimization ToolBox. User Guide},
The MathWorks, Natick, MA, 1999.


\bibitem{CoRelax}
\textsc{W. Collier}:
\textit{Estimating stresses in a  partially inflated high altitude balloon
using~a relaxed energy approach},
\textbf{61} No. 1  17-40 (2003).



\bibitem{dacorogna} \textsc{B. Dacorogna}:
\textit{Direct Methods in the Calculus
of Variations}, Springer-Verlag, New York,  1989.


\bibitem{doCarmo}
\textsc{M. P. do Carmo}: \textit{Differential Geometry of Curves and
Surfaces},
Prentice-Hall,
Englewood Cliffs, NJ, 1976.


\bibitem{evans} \textsc{L. C. Evans}: \textit{Partial Differential Equations},
Graduate Studies in Mathematics, Vol. 19, American Mathematical
Society,  2002.


\bibitem{fisher} \textsc{D. Fisher}:
\textit{\em Configuration dependent pressure potentials},
Journal of Elasticity, \textbf{19} (1988) 77-84.



\bibitem{GiTr} \textsc{D. Gilbarg \& N. S. Trudinger}:
\textit{Elliptic Partial Differential Equations of Second Order},
Springer-Verlag, Berlin, 1983.


%



\bibitem{jonesA} \textsc{W. V. Jones}:
\textit{\em Evolution of Scientific Ballooning},
29th International Cosmic Ray Conference,  Pune \textbf{10} 173-184  (2004).

\bibitem{jonesB} \textsc{W. V. Jones}:
\textit{\em Pioneering space research with balloons},
COSPAR~2006-A-03074,
36th COSPAR Scientific Assembly,
16-23 July 2006, Beijing, China.


\bibitem{libaisimm} \textsc{A. Libai
\& J. G. Simmonds}:
\textit{The nonlinear theory of elastic shells},
2nd Edition,
Cambridge University Press, New York,
1998.



\bibitem{Pipkin}
\textsc{ A. C. Pipkin}:
\textit{ Relaxed energy densities for large deformations of
 membranes}, IMA Journal of Applied Mathematics,
 \textbf{52} 297-308 (1994).




\bibitem{Smc} \textsc{J. H. Smalley}:
\textit{ Development of the e-balloon},
AFCRL-70-0543,
National Center for Atmospheric Research, Boulder, Colorado,
June 1970.


\bibitem{SteigPip} \textsc{ D. J. Steigman \& A. C. Pipkin}:
\textit{ Axisymmetric tension fields}, ZAMP,
\textbf{40}  526-542 (1989).



\bibitem{SteinHedge}
\textsc{M. Stein \& J. M. Hedgepath}:
\textit{Analysis of partially wrinkled membranes},
NASA TN D-813, 1961.

\bibitem{taylor} \textsc{G. I. Taylor}:
\textit{On the stability of parachutes}, in
The Scientific Papers of Sir G. I. Taylor, Vol.~III,
ed. G. K. Batchelor,
Cambridge University Press, 1963.



\end{thebibliography}
\end{document}